%% file: tziolas.tex
\documentclass[12pt,twoside]{gen-j-l}
\usepackage{amsthm}
\usepackage{amsmath}
\usepackage{amstext}
\usepackage{amsfonts}
\usepackage{latexsym}
\usepackage{amscd}
\usepackage{xypic}
\newcommand{\la}{\ensuremath{\longrightarrow}}

\newcommand{\pone}{\ensuremath{\mathbb{P}^{1}}}

\newcommand{\C}{\ensuremath{\mathrm{C}}}

\newcommand{\sheaf}{\ensuremath{\mathcal{O}}}

\newcommand{\ts}{\ensuremath{\textstyle}}
\newtheorem{theorem}{Theorem}[section]
\newtheorem{lemma}[theorem]{Lemma}
\newtheorem{definition}[theorem]{Definition}
\newtheorem{proposition}[theorem]{Proposition}
\newtheorem{corollary}[theorem]{Corollary}

\begin{document}
\title[Divisorial extremal neighborhoods]{Three dimensional divisorial extremal neighborhoods}
\author{Nikolaos Tziolas}
\address{Max Planck Institute for Mathematics, Vivatsgasse 7, Bonn 53111, Germany }
\address{mailing address: Department of Mathematics, University of Crete, Knossos Avenue, Heraklion, 71409, Greece}
\email{tziolas@mpim-bonn.mpg.de}





\begin{abstract}
In this paper we study divisorial extremal neighborhoods $C \subset Y \stackrel{f}{\la} X \ni 0$, 
such that $0 \in X$ is a $cA_n$ type threefold terminal singularity, and $\Gamma =f(E)$ is a smooth curve, 
where $E$ is the $f$-exceptional divisor. We view a divisorial extremal neighborhood as a one parameter 
smoothing of certain surface singularities, and based on this we give a classification of such neighborhoods.
\end{abstract}

\maketitle

\input{introduction}
\input{section0}

\input{section1}

\input{section2}

\input{section3}

\input{section4}

\input{bibliography}
\end{document}

%% file: introduction.tex
\section{Introduction}
One of the most important and difficult problems of modern three dimensional birational geometry, is the study of the structure 
of birational maps between two Fano-Mori fiber spaces $X/S$ and $X^{\prime}/S^{\prime}$. 
Many of the classical rationality problems fall into this context. 

In order to achieve this, the Sarkisov program was developed by Corti, Reid and Sarkisov~\cite{Cor95}, 
whose aim is to factorize any birational map between Fano-Mori fiber spaces as a composition of 
simpler maps, the so called ``elementary links''. 
These links consist of divisorial contractions, flips and flops, and their structure is at the moment 
not well understood. Flops were classified 
by Koll\'{a}r~\cite{Ko91}, and flips by Koll\'{a}r and Mori~\cite{Ko-Mo92}. The structure of divisorial contractions is still an open 
problem. Therefore, and in order to be able to have a working form of the Sarkisov program, it is essential to try and classify 
divisorial contractions $E \subset Y \stackrel{f}{\la} \Gamma \subset X$. 

Mori and Cutkosky~\cite{Cut88a} have classified such contractions in the case that $Y$ is Gorenstein. The case that $\Gamma$ is just a point, 
has been studied by Luo, Corti, Kawakita and others. 

The case that $\Gamma$ is a curve have been studied by Kawamata~\cite{Kaw94}, and Tziolas~\cite{Tzi02}. Kawamata showed that 
if there is a cyclic quotient singularity $P \in \Gamma \subset  X$, then $\Gamma=\{ P\}$, and $f$ is a weighted blow up. 
In~\cite{Tzi02}, the case when $\Gamma$ is a smooth curve and the general section $S$ of $X$ through $\Gamma$
is a $D_{2n}$ singularity was studied, and a classification was given. 

In this paper we will give a classification of divisorial contractions in the case that $\Gamma$ is a smooth curve, 
and $X$ has $cA_n$ type singularities. Hence this paper together with~\cite{Tzi02} and a forthcoming one that 
will complete the study of the $cD$ cases, will give a complete classification of 
divisorial contractions $E\subset Y \stackrel{f}{\la} \Gamma \subset X$, when $X$ is Gorenstein and $\Gamma$ a smooth curve.

The methods used are completely different than the ones used in~\cite{Tzi02}. 
The proper setting of the problem is the local one. We will study maps $C \subset Y \stackrel{f}{\la} X \ni 0$ such that 
$0\in X$ is the germ of a $cA_n$ type singularity, the exceptional set of $f$ is an irreducible divisor and $C=f^{-1}(0)$ an 
irreducible curve. We call such a map a divisorial extremal neighborhood.
 
We view a divisorial extremal neighborhood as a one parameter deformation of certain surface singularities, and then we reduce the problem 
of understanding contractions in understanding deformations of surface singularities. The idea is the following. 
Start with an index 1 threefold singularity $0 \in X$. Then in suitable analytic coordinates, it can be described 
by $(g(x,y,z)+tf(x,y,z,t)=0)\subset \mathbb{C}^4$, where $(g(x,y,z)=0)\subset \mathbb{C}^3$ is a DuVal surface singularity. 
Therefore $X$ can be viewed as a one parameter deformation of a DuVal singularity. Let now $T$ be the general section of $X$ through $0$, 
and $Z=f^{\ast}T$. Then $Y$ is a one parameter deformation of $Z$. Moreover, since $Y$ itself is terminal and $\mathbb{Q}$-Gorenstein 
(i.e., $K_Y$ is $\mathbb{Q}$-Cartier), $Y$ is a $\mathbb{Q}$-Gorenstein, terminal smoothing of the surface singularity $Z$. 

Therefore we may construct a divisorial extremal neighborhood as follows. Start with a birational map 
$f \colon C \subset Z \la 0 \in T$, where $0 \in T$ is DuVal, $C \cong \mathbb{P}^1$, and $-K_Z$ is $f$-ample. 
Let $Y$ be a one parameter $\mathbb{Q}$-Gorenstein, terminal smoothing of $Z$. Then $f$ extends to a birational morphism 
$Y \la X$~\cite{Ma-Ro71}~\cite[Proposition 11.4]{Ko-Mo92}, where $X$ is a one parameter deformation of $T$. 
$X$ is cDV and therefore the contraction is divisorial and not flipping. 
This way we obtain a three dimension divisorial extremal neighborhood. 

The success of this method depends on showing that given a divisorial contraction as above, $Z=f^{\ast}T$ has reasonably good 
singularities whose deformation spaces can be worked and therefore classify all divisorial extremal neighborhoods in this way. 

Let $\Gamma \subset S \subset X$ be the general section of $X$ through $\Gamma$. This must be DuVal~\cite{Ko-Mo92}. 
 Let $C=f^{-1}(0)$, be the central curve, 
and $S_Y=f_{\ast}^{-1}S$. Then we will proceed by considering two cases. The first one is when $C \subset S_Y$. This will be called, 
according to the terminology introduced by Koll\'{a}r and Mori~\cite{Ko-Mo92}, a semistable extremal neighborhood. 
The other case, i.e., when $C \not \subset S_Y$, is the nonsemistable case. In the nonsemistable case, two subcases will have to 
be considered. When the general member $Z$ of $|\sheaf_Y|$ is normal, and when it is not. Theorems~\ref{normal1},~\ref{normal2} are 
treating the first case, and theorems~\ref{nonnormal1},~\ref{nonnormal2}, treat the second case. 
Finally theorems~\ref{semistable1},~\ref{semistable2}, treat the semistable case, and corrollary~\ref{Pell} shows that there 
is a one to one correspondence between semistable extremal neighborhoods and the solutions of a certain Pell equation. As an application 
we consider the cases that $0 \in X$ is a $cA_1$ and $cA_2$ singularities.

Most of this work was made during my stay at Max Planck Institute f\"{u}r Mathematik. I would also like to thank J\'{a}nos Koll\'{a}r for many fruitful discussions during my visit to Princeton University in April 2003, as well as Princeton University for the hospitality.

%% file: section0.tex
\section{Terminology, Notations.}

We start by defining divisorial contractions and extremal neighborhoods.

\begin{definition}
A threefold divisorial contraction is a proper morphism $E \subset Y \stackrel{f}{\la} \Gamma \subset X$, with the following 
properties:
\begin{enumerate}
\item $Y$ is $\mathbb{Q}$-factorial.
\item $Y-E \cong X-\Gamma$, and $E$ is an irreducible divisor.
\item $-K_Y$ is $f$-ample. 
\item $X$ and $Y$ have terminal singularities.
\end{enumerate}
\end{definition}
\begin{definition}
A three dimensional extremal neighborhood is a proper morphism $C \subset Y \stackrel{f}{\la} X \ni 0$, 
with the following properties:
\begin{enumerate}
\item $Q \in X$ is the germ of a terminal singularity, with $Y$ also terminal.
\item $C=f^{-1}(0)$ is an irreducible curve.
\item $-K_Y$ is $f$-ample.
\end{enumerate}
\end{definition}
If the exceptional set of $f$ is an irreducible divisor, then the extremal neighborhood is called divisorial. 
Otherwise it is called flipping. So, a divisorial extremal neighborhood is the local version of a divisorial contraction when 
$\Gamma$ is a curve.

Let $C \subset Y \stackrel{f}{\la} X \ni 0$ be a divisorial extremal neighborhood. We will denote by $E$ the 
$f$-exceptional divisor, and by $\Gamma$ it's center on $X$. i.e., $\Gamma = f(E)$. In this paper $\Gamma$ will 
always be considered to be a smooth curve. 

Let $S$ be the general section of $X$ through $\Gamma$, which must be DuVal~\cite{Ko-Mo92}. Then an extremal neighborhood will 
be called semistable if $C \subset S_Y=f_{\ast}^{-1}S$, and non-semistable if $C \not\subset S_Y=f_{\ast}^{-1}S$.

A threefold $Y$ will be called $\mathbb{Q}$-Gorenstein, iff $K_Y$ is $\mathbb{Q}$-Cartier.
Finally we will define the notion of $T$-singularities.

\begin{definition}[Definition 3.7~\cite{KoBa88}]
A normal surface singularity is called a $T$-singularity, if it is a quotient singularity and admits a $\mathbb{Q}$-Gorenstein 
one parameter smoothing. Such a smoothing must necessary be terminal.
\end{definition}
Such singularities have been completely classified~\cite{Bri68},~\cite{KoBa88}. We will also use the notion of 
semi-log-canonical (\textit{slc}), and semi-log-terminal (\textit{slt}) singularities as they appear in~\cite[Definition 4.17]{KoBa88}.

%% file: section1.tex
\section{Singularities of the general member of $|\sheaf_Y|$.}
Let $C \subset Y \stackrel{f}{\la} X \ni 0$ be a divisorial extremal neighborhood with $0 \in X$ of type $cA_n$, for some $n$, 
and let $T$ be the general section of $X$ through $0$. In this section we will study the singularirities of $Z=f^{\ast}T$.
It turns out that unlike the $cD$ cases where $Z$ has bad singularities, the $cA$ cases behave much better. The next lemma 
shows that.

\begin{lemma}\label{singularities1}
Let $C \subset Y \stackrel{f}{\la} 0 \in X$ be a divisorial extremal neighborhood. Suppose that 
$0 \in X$ is $cA_n$ type, and let $0 \in T \subset X$ be the general section through the singularity. 
Let $Z=f^{\ast}T$. Then $Z$ has semi-log-canonical (slc) singularities. Moreover, the points of index $>1$ are 
semi-log-terminal (slt). 

\end{lemma}

Singularities of this type have been classified by Brieskorn~\cite{Bri68},   
Shepherd-Barron and Koll\'{a}r~\cite{KoBa88}, 
and hence it is possible to describe them precisely. 

\begin{corollary}\label{singularities2}
With assumptions as in the previous lemma, the index $>1$ points of $Z$ are one of the following two types:
\begin{enumerate}
\item Normal klt singularities given as a quotient $W/\mathbb{Z}_n$, where $W$ is $xy-z^{nd}=0$ in $\mathbb{C}^3$, 
the group action is $x \mapsto \zeta x$, $y \mapsto \zeta ^{-1}y$, $z \mapsto \zeta ^a z$, where $\zeta$ is an $n$-th root of unity, 
and $(a,n)=1$.
\item $(xy=0)\subset \mathbb{C}^3 / \mathbb{Z}_n$, and the group action is $x \mapsto \zeta^a x$, $y \mapsto \zeta ^{-a}y$, 
$z \mapsto \zeta z$, where $\zeta$ is an $n$-th root of unity, and $(a,n)=1$.
\end{enumerate}
The index 1 points are either DuVal, degenerate cusps, normal crossing or pinch points.
\end{corollary}
In individual cases we will be able to get more precise information.

\begin{proof}
Let $E$ be the exceptional divisor and $\Gamma = f(E)$. Let $\Gamma \subset S \subset X$ be the general section 
of $X$ containing $\Gamma$. Then by~\cite{Ko-Mo92} this is DuVal, and by~\cite{Tzi02} it must in fact be of type $A_m$. 
Then if $T$ is the general section of $X$ through $0$, the pair $(X,S+T)$ must be log-canonical. This follows from inversion of adjunction. 
$(X,S+T)$ is log-canonical iff $(S,T|S)$ is log-canonical. But if $S$ is given by $xy-z^m=0$, then $T|S$ is just $z=0$, 
and therefore the union of the lines $l_1: x=z=0$, and $l_2:y=z=0$. The pair $(S, l_1+l_2)$ is now easily checked to be log-canonical. 

Now $f^{\ast}S=S_Y+E$, and $K_Y=f^{\ast}K_X+E$. Therefore, \[
K_Y+S_Y+Z=f^{\ast}(K_X+S+T), \]
and therefore $(Y,S_Y+Z)$ is also log-canonical. It is also true that $K_Y+S_Y=f^{\ast}(K_X+S)$. From this it follows that 
the high index points of $Y$ must be on $S_Y$, and it immediately follows again from inversion of adjunction that the high index points of 
$Z$ must be semi-log-terminal (slt). If it happens that $Z$ is normal, then they are just klt quotient singularities. 
However, it is possible that $Z$ is not normal and therefore it can be just slt, as we will see later.

Now if it happens that $Z$ is normal, then the high index points are klt singularities that have one parameter 
$\mathbb{Q}$-Gorenstein smoothings. These are classified~\cite[Proposition 3.10]{KoBa88} as mentioned above and are those given 
in part $(1)$ of the corollary. 

The index 1 points are log-canonical. Moreover, it is easy to see that $R^1f_{\ast}\sheaf_{Z} =0$, and hence they are also rational 
singularities. But then by~\cite[lemma 9.3]{Kaw88}, they must be DuVal.

Suppose now that $Z$ is not normal. Once again the high index points are slt with one parameter 
$\mathbb{Q}$-Gorenstein smoothings. They are classified~\cite[Theorem 4.23, 5.1]{KoBa88}, and are exactly what stated in the corollary.
\end{proof}
The most general result that we can give towards classification of divisorial extremal neighborhoods $Y \la X$, with $X$ $cA_n$ is the following.
\begin{theorem}\label{deformation}
Fix an $A_n$ surface germ $(0 \in T)$, and the germ $(C \subset \Sigma)$ of a surface $\Sigma$ along a smooth curve $C$. 
Then there exist a closed subscheme $W_{C,\Sigma} \subset Def(0 \in T)$ with the following property. 
Let $0 \in \Gamma \subset X$ be a threefold germ along a smooth curve $\Gamma$, such that the general section of $X$ at $0$ is 
isomorphic to $(0 \in T)$, and $(\Gamma \subset S) \cong (\C \subset \Sigma)$, where $S$ is the general section of $X$ through $\Gamma$. 
$X$ is a deformation of $T$ over the unit disk $\Delta$, and there is a natural map $\Delta \la Def(0 \in T)$. 

Then a divisorial extremal neighborhood $Y \la X$ contracting an irreducible divisor $E$ onto $\Gamma$ exists, iff
\[
\text{Im}[\Delta \la Def(0\in T)] \subset W_{C,\Sigma}.\]
\end{theorem}

The proof of the theorem will be given at the end of section 5.

%% file: section2.tex
\section{Non-semistable neigborhoods.}

In this section we will classify the non-semistable divisorial extremal neighborhoods $f \colon C \subset Y \la 0 \in X$.
Let $E$ be the exceptional divisor, and $\Gamma =f(E)$. Let $S$ be the general section of $X$ through $\Gamma$. Then $C \not \subset S_Y$. 
Let $0 \in T$ be the general section of $X$ through $0$, which is an $A_m$ singularity for some $m$, and let $Z=f^{\ast}T$. 
The following two examples will show that there are cases when $Z$ is normal, and others that it is not. Therefore we will have to treat 
the two cases differently. 

\textbf{Example 1.} Start with the germ of a smooth surface along four $\mathbb{P}^1$'s with the following configuration: \[
\begin{array}{c}
\stackrel{-3}{\circ} \mbox{\noindent ---} \stackrel{-2}{\circ} \mbox{\noindent ---} \stackrel{-3}{\circ} \\
| \\
\bullet \\
\mbox{\small{-1}}
\end{array} \]
Now contract all except the $-1$ curve to get a birational map $ C \subset Z \stackrel{f}{\la} T \ni 0$, where $0 \in T$ is an 
$A_2$ DuVal singularity, and $Z$ has exactly one singular point $P$, which is a cyclic quotient singularity of type $1/12 (1,5)$. 
This is $\frac{1}{n^2d}(1,and-1)$, with $n=2$, $d=3$, $a=1$, and hence it is a $T$-singularity\cite{KoBa88} and hence admits a 
$\mathbb{Q}$-Gorenstein terminal smoothing $Y$, which must be terminal by~\cite[Corollary 3.6]{KoBa88}. Now the map $f$ extends to a 
birational morphism $Y \la X$, which is a divisorial extremal neighborhood such that $Z=f^{\ast}T$ is normal.

The following example which was communicated to me by Shigefumi Mori, shows that it is indeed possible to have $Z$ non-normal.

\textbf{Example 2.} Consider the configuration \[
\stackrel{-2}{\circ} \mbox{\noindent ---} \stackrel{-4}{\circ} \mbox{\noindent ---} 
\stackrel{-4}{\bullet} \mbox{\noindent ---} \stackrel{-1}{\bullet} \mbox{\noindent ---} 
\stackrel{-2}{\circ} \mbox{\noindent ---} \stackrel{-2}{\circ} \mbox{\noindent ---} \stackrel{-3}{\circ} 
\]
Contracting all curves except those marked by a solid circle, we obtain a morphism $\overline{Z} \stackrel{\overline{f}}{\la} T$, where 
$T$ is an $A_2$ singularity, the exceptional set is the union of two rational curves $C_1$, $C_2$, and 
$\overline{Z}$ has exactly two singular 
points $P_1$ and $P_2$. From the construction it is clear that $(P_1 \in \overline{Z}) \cong 1/7(1,2)$, 
and $(P_2 \in \overline{Z}) \cong 1/7(1,-2)$. 
The curves $C_1$ and $C_2$ can be identified to a rational curve $C$, and we see $\overline{Z}$ as the normalization of $Z$ which 
has exactly one high index point which is of type $(xy=0)/\mathbb{Z}_7(2,-2,1)$. 

$Z$ has a $\mathbb{Q}$-Gorenstein terminal smoothing $Y$, and the map $Z\la T$ extends to a map $Y \stackrel{f}{\la}X$, that is 
a divisorial extremal neighborhood such that $Z=f^{\ast}T$ is not normal. 

Hence from the above examples it is clear that we must consider separately the cases of non-semistable neighborhoods 
with $Z$ normal and non-normal. We will call the first kind \textit{Normal non-semistable neighborhoods}, and the second 
\textit{Non-normal non-semistable neighborhoods}

\subsection{Normal non-semistable neighborhoods}

We now want to classify non-semistable divisorial neighborhoods $C \subset Y \stackrel{f}{\la} X \ni 0$,  
such that if $T$ is the general section through $0$, then $Z=f^{\ast}T$ is normal. 
Let $S$ the general section of $X$ through $\Gamma$. Since we are in the non-semistable case, $S_Y \cong S$.
First we want to get more information about the singularities of $Z$. 

\begin{lemma}
With assumptions as above, $Z$ has exactly one singular point $P$ of index bigger than 1, and at most one index 1 singular point which 
must be DuVal. Moreover, $(P \in Z) \cong \frac{1}{n^2d} (1,and-1)$, with $(a,n)=1$.
\end{lemma} 

\begin{proof}
By lemma~\ref{singularities1} and corollary~\ref{singularities2}, it follows that the high index points are of the type stated in the lemma, 
and they must be on $S_Y \cap C$, which is a single reduced point from the assumption that $S_Y\cong S$. Therefore there is 
only one high index point.

The index 1 singular points must be DuVal by lemma~\ref{singularities2}. The only thing left is to show that there can be at most one. 
Let $U \stackrel{g}{\la} Z$ be the minimal resolution of $Z$. Let $C_U=g^{-1}_{\ast}C$. Since $0 \in T$ is $DuVal$, $U$ must dominate 
the minimal resolution of $0 \in T$, and therefore it is obtained from it by a sequence of blow ups. Therefore, $C_U^2=-1$. 
If there are more than two DuVal points on $Z$, then in the dual graph of $Z$ there will be the configuration \[
\stackrel{-2}{\circ} \mbox{\noindent ---} \stackrel{-1}{\circ} \mbox{\noindent ---} \stackrel{-2}{\circ} \]
After contracting $C_U$ we get the configuration $\stackrel{-1}{\circ} \mbox{\noindent ---} \stackrel{-1}{\circ} $, which is not contractible 
anymore. Therefore there must be at most one DuVal point on $Z$.
\end{proof}

Hence we must classify surface birational maps $C \subset Z \stackrel{f}{\la} T \ni 0$, such that $0 \in T$ is $A_m$, for some $m$, 
and $Z$ has one $T$-singularity as above and at most one DuVal point. Moreover, we want that the resulting threefold extremal neighborhood 
$Y \la X$ is divisorial, and if $E$ is the exceptional divisor, then $\Gamma =f(E)$ is smooth. 

First we show how to compute the multiplicity of $\Gamma$. 

\begin{lemma}\label{multiplicity}
Let $C \subset Z \stackrel{f}{\la} T \ni 0$ be a birational map of surfaces. Let $Y$ be a $\mathbb{Q}$-Gorenstein terminal smoothing 
of $Z$, and let $Y \stackrel{f}{\la} X$ be the threefold extremal neighborhood that is produced by extending $f$ in the family. Assume that 
it is divisorial, $E$ the exceptional divisor, and $\Gamma =f(E)$. Then the multiplicity $m_{\Gamma}$ of $\Gamma$ is given by \[
m_{\Gamma}=-\frac{(C \cdot K_Z)^2}{C^2} \]
\end{lemma}

\begin{proof}
Look at how the curve $C$ degenerates in the family. In the general fiber $Y_g$, the exceptional set is a disjoint union of $m=m_{\Gamma}$ 
smooth rational $-1$-curves, say $C_1, \cdots , C_m$. They degenerate to $\nu C$ in the central fiber $Z$. 
Then we have 
\begin{eqnarray*}
\nu C \cdot K_Z = (\sum_{i=1}^m C_i)\cdot K_{Z_g}= -m \\
\nu^2 C^2= (\sum_{i=1}^m C_i)^2= -m
\end{eqnarray*}

Now suppose that $K_Z=f^{\ast}K_T+aC$. This implies that \[
a=\frac{C\cdot K_Z}{ C^2}=\nu.\]
 Therefore we obtain that \[
m_{\Gamma}=-\frac{(C \cdot K_Z)^2}{C^2} \]
as claimed.
\end{proof}
We are now in position to describe all normal non-semistable divisorial extremal neighborhoods. 
\begin{theorem}\label{normal1}
Normal non-semistable divisorial extremal neighborhoods $C \subset Y \stackrel{f}{\la} X \ni 0$, 
with $X$ of type $cA_m$ for some $m$, and $\Gamma$ smooth,
are in one to one correspondence 
with $\mathbb{Q}$-Gorenstein smoothings of normal surface germs $C \subset Z$ that satisfy the following conditions:
\begin{enumerate}
\item $Z$ has exactly one high index point $P$ and at most one index one singular point. Moreover, $(P\in Z) \cong \frac{1}{n^2d}(1,and-1)$, 
with $(a,n)=1$, and the index 1 point is DuVal.
\item $K_Z \cdot C =-1/n$, and $C^2=-1/n^2$.
\end{enumerate}
As a consequence of the second condition follows that $K_Y\cdot C=-1/n$, and $n=\text{index}(Y)$.

\end{theorem}

To complete the description of the resulting neighborhood 
 $Y\stackrel{f}{\la} X$, we want to obtain information about the singularity of $X$, and of the general 
section $S$ of $X$ containing $\Gamma$, from data of $Z$.

\begin{theorem}\label{normal2}
Let $C \subset Z$ be a normal surface germ with the properties stated in theorem~\ref{normal1}, and let $C \subset Z \stackrel{f}{\la} T \ni 0$ 
be the contraction of $C$. Let $C \subset Y \stackrel{f}{\la} X \ni 0$ be the associating threefold contraction. Then
\begin{enumerate}
\item $C \subset Y \stackrel{f}{\la} X \ni 0$ can be obtained from the following diagram
\[
\xymatrix{ 
            & V \ar[dl]_{\nu} \ar[dr]^{\pi}\\
W \ar[dr]_g &   & Y \ar[dl]^f \\
            & X & }
\]
where $W$ is the blow up of $X$ along $\Gamma$. There are two $g$-exceptional divisors, a ruled surface $E$ over $\Gamma$, 
and $F\cong \mathbb{P}^2$ over $0$. $V$ is the $\mathbb{Q}$-factorialization of $E$, and $\pi$ contracts $F_V \cong \mathbb{P}^2$ to a point.
\item $0 \in X$ is $cA_{d-1}$, where $d$ is the torsion of $Cl^{sc}(Z)/([C]=0)$ and $Cl^{sc}(Z)$ is the group of $\mathbb{Q}$-Cartier divisors 
of $Z$.
\item Let $S$ be the general section of $X$ containing $\Gamma$. It is $A_{\nu}$ for some $\nu$. Suppose that $\Gamma$ intersects the 
$E_k$ exceptional curve in the fundamental cycle of the minimal resolution of $S$. Then,
\[
n=\frac{\nu +1}{(k,\nu +1)}
\]
\end{enumerate}
\end{theorem}
\begin{proof}
Part $1.$ of theorem~\ref{normal1} follows from lemma~\ref{singularities1}. 
For the second part we must find the conditions for $0 \in T$ to be $DuVal$, 
and $\Gamma$ smooth. From lemma~\ref{multiplicity}, it follows that $\Gamma$ is smooth iff 
\begin{equation}
C^2+(K_Z\cdot C)^2=0. 
\end{equation}
We now need the following characterization of cDV points.
\begin{lemma}
Let $C \subset Y \stackrel{f}{\la} X \ni 0$ be a divisorial extremal neighborhood as above. Then $X$
has index 1, iff one of the following 
equivalent conditions is satisfied.
\begin{enumerate}
\item $Cl^{sc}(Y)$ is torsion free.
\item $K_Y\cdot C=-1/n$, where $n$ is the index of $Y$.
\end{enumerate}
\end{lemma}
\begin{proof}
Suppose that $Cl^{sc}(Y)$ is torsion free. Then $f^{\ast}K_X=K_Y-E$, where $E$ is the $f$-exceptional divisor. Let $m$ be the index of $X$. 
Then $mK_X=0$, and hence $m(K_Y-E)=0$. Therefore $K_Y-E=0$, and hence $m=1$. 

Conversely, suppose that $X$ has index 1 and hence it is cDV. Let $D_Y \in Cl^{sc}(Y)$ be a torsion element, i.e., $mD_Y=0$, 
for some $m$, and 
hence $D_Y \cdot C=0$. Then $D_Y = f^{\ast} D_X$ for some $\mathbb{Q}$-Cartier divisor $D_X$ in $X$. But since $X$ is cDV, $D_X$ must in fact 
be Cartier~\cite{Kaw88}, and hence $D_X=0$ and therefore $D_Y=0$ as well. Hence $Cl^{sc}(Y)$ is torsion free.

To see the second assertion now. By~\cite[Corollary 1.10]{Mo88}, $Cl^{sc}(Y)$ is torsion free iff there exist an effective divisor $D \in Cl^{sc}(Y)$ 
such that $D\cdot C=1/n$. On the other hand, from the proof of lemma~\ref{multiplicity} it follows that $K_Y \cdot C =E\cdot C=-1/\nu$, for 
some number $\nu | n$. Then $(nD+\nu E)\cdot C =0$, and hence $nD+\nu E = f^{\ast} (nD_X)$ for some divisor $D_X$. But since $f$ is generically 
the blow up of $\Gamma$, it follows that $\nu \geq n$, and hence $\nu = n$. 
\end{proof}
Therefore, $K_Y \cdot C = -1/n$ and hence $(1)$ is also equivalent to the condition $C^2=-1/n^2$. This concludes the proof of theorem~\ref{normal1}. 

To see part $(1)$ of theorem~\ref{normal2} now. By~\cite{Tzi02}, it follows that the contraction $Y\la X$ can be constructed by the following 
diagram:
\[
\xymatrix{ 
V \ar[d]_{\nu} \ar@{-->}[rr]^{\phi} & & V^{\prime} \ar[d]^{\pi}\\
W \ar[dr]_g &   & Y \ar[dl]^f \\
            & X & }
\]
where, $W$ and $V$ are as in the statement of the theorem, and $\phi$ is a composition of flips. We want to show that in this case there 
are no flips and in fact $V=V^{\prime}$. Suppose that there are flips, and let $C^{\prime}$ be the last flipped curve. Then 
$K_{V^{\prime}}\cdot C^{\prime} > 0$. Let $h=f \circ \pi$. Then \[
K_{V^{\prime}}+S_{V^{\prime}}=h^{\ast}(K_X+S) \]
where $S$ is the general section of $X$ through $\Gamma$. Then it follows that $S_{V^{\prime}} \cdot C^{\prime} <0$, and 
therefore $C^{\prime} \subset S_{V^{\prime}}$. $K_{V^{\prime}}$ is $\pi$-negative and 
hence $C^{\prime}$ does not contract, and $\pi (C^{\prime})=C$. 
But this implies that $C \subset S_Y$, which is not the case since we only consider non-semistable neighborhoods. Therefore there are no flips 
and hence $V=V^{\prime}$.

We now want to show that $F_V$ is contracted to a point by $\pi$. By~\cite{Tzi02} it follows that $K_W+S_W=g^{\ast}(K_X+S)$. 
Therefore $K_V+S_V=h^{\ast}(K_X+S)$. Moreover, $K_Y+S_Y=f^{\ast}(K_X+S)$. Hence we see that \[
K_V+S_V=\pi^{\ast}(K_Y+S_Y).\]
Therefore $F_V$ is $\pi$ crepant. If $C \not \subset S_Y$, and $F_V$ contracts onto $C$, then in fact $F_V$ is crepant over $Y$. 
But this is impossible since $Y$ is assumed to be terminal. Hence, and since $C \not \subset S_Y$, $F_V$ must contract to a point. 
In particular this implies that no $\nu$-exceptional curve $\Delta$ can be contained in $F_V$. If they did then they would have 
to be contracted by $\pi$. But this is impossible since $K_V$ is $\pi$-negative and $K_V \cdot \Delta =0$. 
Hence $F_V \cong F \cong \mathbb{P}^2$. 

Now look at $f \colon Z \la T$. If the conditions of theorem~\ref{normal1} are satisfied, then $0 \in T$ is DuVal, and in particular of 
type $A_{d-1}$. We want to find $d$. Suppose that $0 \in T$ is given by $xy-z^d=0$, and let $l \subset X$ be the line $x=z=0$. 
Then $d=\text{index}(l)$. Let $l^{\prime}=f^{-1}_{\ast}l$. Then $f^{\ast}l=l^{\prime}+\delta C$, and since $dl=0$, it follows that 
$d(l^{\prime}+\delta C)=0$, or equivalently that $dl^{\prime}+mC=0$, for $d,m \in \mathbb{N}$. 
Now it is known that $Cl^{sc}(Z)$ is generated by $[l]$ and $[C]$, with the relation $dl^{\prime}+mC=0$. 
Therefore it follows that $d$ is just the torsion part of $Cl^{sc}(Z)/([C]=0)$. 

Finally we want to get information about the general section $S$ of $X$ through $\Gamma$. It is of type $A_{\nu}$, for some $\nu$~\cite{Tzi02}. 
Let $E$, $F\cong \mathbb{P}^2$ be the $g$-exceptional divisors, and $E_V$, $F_V$ their birational transform in $V$. 
Let $l \subset F$ be a general line. Then

\textbf{Claim:}

\begin{enumerate}
\item \[
\begin{array}{ccc}
\text{index}(E_V)=\frac{\nu -k +1}{(k,\nu +1)} & , & \text{index}(F_V)=\text{lcm}(k,\nu -k+1)
\end{array}
\]
\item \[
l \cdot F_V = -\frac{\nu +1}{k(\nu -k+1)}
\]
\end{enumerate}
Now suppose that $K_V=\pi^{\ast} K_Y + a F_V$. $K_V \cdot l = K_W \cdot l = -1$. Hence $a = k(\nu -k+1)/(\nu +1)$, and therefore \[
K_V=\pi^{\ast} K_Y +\frac{k(\nu -k+1)}{\nu +1}F_V.\]
It now immediately follows that \[
n=\text{index}(Y)=\frac{\nu +1}{(k,\nu +1)} \]
as claimed. 

\textbf{Proof of the claim.} We start by computing the indices of $E_V$ and $F_V$. The following result will be very useful.

\begin{lemma}
Let $0 \in S$ be an $A_n$ DuVal singularity. Let $U \stackrel{f}{\la} S$ be the minimal resolution and $E_i$, $i=1, \ldots n$, the exceptional 
curves. Let $l$ be a line in $S$ that its birational transform $l^{\prime}=f_{\ast}^{-1}l$ intersects $E_k$. Then the index of $l$ in $S$ 
is $(n+1)/(k,n+1)$.
\end{lemma}

\textbf{Remark:} The above lemma has the following interpretation. Suppose that $0 \in S$ is given by $xy-z^{n+1}=0$. Let $l$ be the line 
given by $x-z^k=y-z^{n-k+1}=0$. Then the index of $l$ is $(n+1)/(k,n+1)$. Under suitable coordinates, all smooth curves through the singularity 
are given by equations of the above form~\cite{Jaf92}, and hence we have described the index of any smooth curve depending on its position in 
the fundamental cycle in the minimal resolution.

\begin{proof}[Proof of the lemma.]
It is not difficult to find that \[
f^{\ast}l=l^{\prime}+\frac{n-k+1}{n+1} \sum_{i=1}^{k} i E_i + \frac{k}{n+1} \sum_{i=1}^{n-k} (n+1-k-i) E_{k+i}. \]
The lemma now follows immediately.
\end{proof}
Look at $W$ now. $L=E \cap F$ is a line in $F=\mathbb{P}^2$~\cite{Tzi02}. At the generic point of $L$, $\text{Spec}(\sheaf_{W,L})$ is 
two DuVal singularities, $E$, $F$ correspond to lines, and the extended dual graphs are by~\cite[Propositions 4.5, 4.6]{Tzi02}
\[
\begin{array}{c}
\stackrel{E_1}{\circ} \mbox{\noindent ---} \stackrel{E_2}{\circ} \mbox{\noindent ---} \cdots \mbox{\noindent ---} 
\stackrel{E_k}{\circ} \mbox{\noindent ---} \cdots \mbox{\noindent ---}\stackrel{E_{n-k}}{\circ} 
\mbox{\noindent ---} \stackrel{F}{\bullet} \\
| \\
\underset{E}{\bullet}
\end{array}
\]
and
\[
\stackrel{F}{\bullet} \mbox{\noindent ---}  \stackrel{E_{n-k+2}}{\circ} \mbox{\noindent ---} 
{\circ} \mbox{\noindent ---} \cdots \mbox{\noindent ---}\stackrel{E_n}{\circ} 
\]
Therefore, $W$ is singular along two lines. $L =E \cap F$, and another that does not lie in $E$ but only in $F$. 
Therefore, by using the previous lemma, we see that at the generic point of $L$, $E$ has index $(\nu -k+1)/(k,\nu +1)$, and $F$ $lcm(k,\nu -k+1)$. 
$V$ is just the blow up of $W$ along $mE$, where $m=(\nu -k+1)/(k,\nu +1)$, and part $(1)$ of the claim follows immediately. 

Now to compute $l \cdot E_V$. for a general line $l \subset F_V \cong \mathbb{P}^2$. $E_V$ has index $(\nu -k +1)/(k,\nu +1)$. 
Hence $ (\nu -k +1)E_V \cdot l = \text{length}(\nu -k +1)E_V \cap l$. To find this we can work at the generic point of $L$, 
where $E$ and $F$ are lines, and calculate $(\nu -k +1)E_V \cap F_V$, or $E_V \cap (\nu -k +1)F_V$. 

By~\cite{Jaf92}, under suitable coordinates, $V$, $E$ and $F$ are given at the generic point of $L$ by $xy-z^{\nu -k +1}=0$, 
$x-z^k=y-z^{\nu -2k +1}=0$ and $x=z^{\nu -k+1}=0$. It now follows that \[
E_V \cap (\nu -k +1)F_V= \text{length} \frac{\mathbb{C}[[x,y,z]]}{(x,z^k,z^{\nu -k+1},y-z^{\nu -2k+1})}=k.\]
Hence \[
l \cdot E_V = \frac{k}{\nu -k+1}.\]
By~\cite[Proposition 4.6]{Tzi02} it follows that $g^{-1}(\Gamma)=E+kF$. Therefore, 
$l \cdot (E_V+kF_V)=l \cdot (E+kF) =-1$. since $W$ is the blow up of $\Gamma$. Therefore,
\[
l \cdot F_V = -1-\frac{k}{\nu -k+1}=-\frac{\nu +1}{k(\nu -k+1)},\]
and the claim and hence the theorem is proved.
\end{proof}
All the quantities that appear in theorem~\ref{normal1} can be computed with respect to $n$, $a$, $d$ and the position of $C$ 
in the fundamental cycle of $P \in Z$, and hence the equations $K_Z \cdot C = -1/n$ and $C^2=-1/n^2$ become equations involving 
 the mentioned invariants of the singularity $P\in Z$. Therefore, divisorial contractions of this type are in 
one to one correspondence with the solutions of two diophantine equations. Unfortunately, $C$ can be anywhere in the fundamental cycle 
of $P \in Z$, which makes a reasonable description of $K_Z \cdot C = -1/n$ and $C^2=-1/n^2$ in terms of $n,\; a,\; d$ difficult to get. 
However, we will see that the semistable case behaves much better in this respect.

The following examples show that the curve $C$ can indeed be anywhere in the fundamental cycle.

\textbf{Example 3.} Consider the configuration of $\mathbb{P}^1$'s \[
\stackrel{-2}{\circ} \mbox{\noindent ---}\stackrel{-4}{\circ} \mbox{\noindent ---}\stackrel{-2}{\circ} \mbox{\noindent ---} 
\stackrel{-2}{\underset{\underset{\underset{ -1}{\bullet}}{|}}{\circ}} \mbox{\noindent ---}
\stackrel{-3}{\circ} \mbox{\noindent ---}\stackrel{-3}{\circ} 
\]
Contract all curves except the $-1$ to get a map $C \subset Z \stackrel{f}{\la} T \ni 0$, such that $0 \in T$ is an $A_3$ DuVal singularity 
and $Z$ has exactly one singular point $P$, which is of type $1/n^2d(1,and-1)$, with $n=5$, $d=4$ and $a=3$. 

Now Let $Y \stackrel{f}{\la} X$ be the corresponding threefold contraction, let $E$ the exceptional divisor and $\Gamma =f(E)$. We want to find the 
multiplicity of $E$. 

Let $U \stackrel{g}{\la}Z$ be the minimal resolution, and $E_i$, $i=1 \ldots 6$, the exceptional curves. By construction $C$ intersects $E_4$. 
It is easy to see that \[
K_U=g^{\ast}K_Z -2/5 E_1 -4/5 E_2 - 4/5 E_3 - 4/5 E_4 - 4/5 E_5 - 3/5 E_6 \]
and hence $K_Z \cdot C =-1/5$. It is also easy to compute that \[
g^{\ast}C=C_U+2/25 E_1+4/25E_2+14/25E_3+24/25E_4+9/25E_5+3/25E_6 \]
and hence $C^2=-1/25$. It now follows from lemma~\ref{multiplicity} that $\Gamma$ is smooth. 

\textbf{Example 4.} Consider the configuration of $\mathbb{P}^1$'s \[
\stackrel{-2}{\circ} \mbox{\noindent ---}\stackrel{-3}{\circ} \mbox{\noindent ---}\stackrel{-2}{\circ} \mbox{\noindent ---} 
\stackrel{-2}{\underset{\underset{\underset{ -1}{\bullet}}{|}}{\circ}} \mbox{\noindent ---}
\stackrel{-4}{\circ} 
\]
As in the previous example we obtain a map $C \subset Z \stackrel{f}{\la} T \ni 0$, such that $0 \in T$ is $A_3$, 
and $P \in Z \cong 1/n^2d(1,and-1)$, with $n=6,\; a=1,\; d=1$. Moreover, in the resulting threefold contraction $Y \stackrel{f}{\la} X$, 
$\Gamma$ is also smooth.

These examples together with example 1 earlier show that indeed there is no restriction on the position of $C$ in the fundamental cycle 
of $P\in Z$, even if the singularity of $0 \in T$ is fixed.

%% file: section3.tex
\subsection{Non Normal, Non-semistable neighborhoods.}

In this section we will classify non-semistable divisorial neighborhoods $C \subset Y \stackrel{f}{\la} X \ni 0$, 
with the property that if $0 \in T$ is the general section of $X$ through $0$, then $f^{\ast}T=Z$ is not normal. 
Such neighborhoods do exist as shown by Example 2. 

\begin{theorem}\label{nonnormal1}
Non normal, non semistable divisorial extremal neighborhoods $C \subset Y \stackrel{f}{\la} X \ni 0$, 
with $0 \in X$ of type $cA_k$ for some $k$, are in one 
to one correspondence with $\mathbb{Q}$-Gorenstein smoothings of non-normal surface germs $ C \subset Z$, whose 
normalization $\overline{Z} \stackrel{\pi}{\la} Z$ has the following properties.
\begin{enumerate}
\item $\pi^{-1}(C)=C_1 \cup C_2$, where $C_1, \; C_2$ are smooth rational curves intersecting transversally. ie., $C$ 
breaks in the normalization of $Z$.
\item $\overline{Z}$ has exactly three singular points $P_1, \; P_2,\; Q$. $P_1$ is on $C_1$ but not on $C_2$, 
$P_2$ is on $C_2$ but not on $C_1$, and $Q \in C_1 \cap C_2$. Moreover,
\begin{enumerate}
\item \[
\begin{array}{cc}
(P_1 \in \overline{Z}) \cong \frac{\ts 1}{\ts n}(1,a) & (P_2 \in \overline{Z}) \cong \frac{\ts 1}{\ts n}(1,-a) 
\end{array} \]
for $n, \; a$, with $(a,n)=1$, and $P \in Z$ is a $slt$ singulatity of type $(xy=0)/\mathbb{Z}_n(a,-a,1)$, where $P=\pi(P_1)
=\pi(P_2)$.
\item \[
(Q \in \overline{Z}) \cong \left\{ 
\begin{array}{cc}
\text{smooth} &  \\
\text{DuVal} &  \text{if} \; m_R(Z)=2 \\
\frac{\ts 1}{\ts m}(1,-b) & b | m+1, \; b \neq 1 \ts{mod} m, \; m_R(Z)=3 \\
\frac{\ts 1}{\ts m}(1,-b) & b | 2m+1, \; b \not \equiv 1\text{mod} m, \; m_R(Z)=4 \\
\frac{\ts 1}{\ts m}(1,-b) & b |m+s+1, \; \text{for some } s\geq 1, \\
                          & \; b\geq 2s+1, \; s+1|b-s, \; \text{and} \; m_Q(Z)=4
\end{array} \right. \]
and $R=\pi (Q) \in Z $ is a degenerate cusp~\cite{Ba83}.
\item The extended dual graph of the resolution of singularities of $\overline{Z}$ is a chain. i.e., it is of the form 
\[
\circ \mbox{\noindent---} \cdots \mbox{\noindent---}\circ \mbox{\noindent---}\bullet  \mbox{\noindent---} 
\circ \mbox{\noindent---} \cdots \mbox{\noindent---}\circ \mbox{\noindent---}\bullet \mbox{\noindent---}
\circ \mbox{\noindent---} \cdots \mbox{\noindent---}\circ 
\]
\end{enumerate}
\item \[
(C_1-C_2)^2=n^2[(C_1 \cdot C_2)^2-C_1^2 C_2^2] \]
\item \[
am-(b+1)n \; | \; (m+n^2,b+1+an) \]
\end{enumerate}
\end{theorem}
\textbf{Remark.} Conditions 3. and 4. of the theorem guarantee that in the resulting threefold contraction $Y \la X$, 
$X$ is cDV, and in fact $cA_k$ for some $k$,  and $\Gamma$ is smooth.

Now let $Y \la X$ be the threefold contraction associated to a surface germ $Z$ as in the previous theorem. 
As in the normal case, we want to obtain information about the singularities of $X$ and the general section $S$ 
of $X$ containing $\Gamma$.
\begin{theorem}\label{nonnormal2}
Let $C \subset Z$ be a surface germ with the properties stated in theorem~\ref{nonnormal1}, and let 
$C \subset Z \stackrel{f}{\la} T \ni 0$ 
be the contraction of $C$. Let $C \subset Y \stackrel{f}{\la} X \ni 0$ be the associating threefold contraction. Then
\begin{enumerate}
\item $C \subset Y \stackrel{f}{\la} X \ni 0$ can be obtained from the following diagram
\[
\xymatrix{ 
            & V \ar[dl]_{\nu} \ar[dr]^{\pi}\\
W \ar[dr]_g &   & Y \ar[dl]^f \\
            & X & }
\]
where $W$ is the blow up of $X$ along $\Gamma$. There are two $g$-exceptional divisors, a ruled surface $E$ over $\Gamma$, 
and $F\cong \mathbb{P}^2$ over $0$. $V$ is the $\mathbb{Q}$-factorialization of $E$, and $\pi$ contracts 
$F_V \cong \mathbb{P}^2$ to a point.
\item $K_Y \cdot C = -1/n$, and $n = \text{index}(Y)$.
\item $0 \in X$ is $cA_{d-1}$, where \[
d=\frac{\ts \nu m}{\ts (\nu m, k)} \]
where $\nu=(C_1-C_2)^2 k$, and $k = nm/(m,b+1)$.
\item Let $S$ be the general section of $X$ containing $\Gamma$. It is $A_{\nu}$ for some $\nu$. 
Suppose that $\Gamma$ intersects the 
$E_k$ exceptional curve in the fundamental cycle of the minimal resolution of $S$. Then,
\[
n=\frac{\nu +1}{(k,\nu +1)}
\]
\end{enumerate}
\end{theorem}

\begin{proof}
Let $C \subset Y \stackrel{f}{\la} X \ni 0$ be a non normal, non-semistable neighborhood. Let $0 \in T$ be the general 
section of $X$ through $0$, and $Z=f^{\ast}T$, which by assumption is not normal. By lemma~\ref{singularities1}, there is exactly 
one high index point $P\in Z$, which must be \textit{slt}, and $Z-P$ is \textit{slc}. Let $\overline{Z} \stackrel{\pi}{\la} Z$, 
be the normalization of $Z$. We want to show that $\pi^{-1}(C)$ breaks to two curves. Suppose that it doesn't, and that 
$\overline{C}=\pi^{-1}(C)$ is irreducible. Then we will show that this cannot happen unless $Z$ itself is normal. 

First we will describe the singularities of $\overline{Z}$. By~\cite[Proposition 4.27]{KoBa88}, $\pi^{-1}(P)=\{P_1,\; P_2\}$, 
with $(P_1 \in \overline{Z}) \cong 1/n(1,a)$, and $(P_2 \in \overline{Z}) \cong 1/n(1,-a)$. 

\textbf{Claim:} $\overline{Z}-\{P_1,\; P_2\}$ is smooth. 

By corollary~\ref{singularities2}, the index one points of $\overline{Z}$ must be either degenerate cusps, or normal crossing, or pinch points. 
The only ones with not smooth normalization are the degenerate cusps. The proof of the claim will come from a detailed study of the possible 
degenerate cusps on $Z$. At first we remark that the normalization of a degenerate cusp is cyclic quotient and hence rational~\cite{Ba83}. 
Therefore $\overline{Z}$ has rational singularities. 

$\overline{Z}$ fits into the following diagram:
\[
\xymatrix{ 
            & \overline{Y} \ar[dl]_{h} \ar[rr]^{p} &      &   Y \ar[dr]^{g}         &    \\
\overline{Z}\ar[drr]_{\overline{f}} &                      &      &             &  Z\ar[dll]^{f}  \\
            &                      &  X   &             & }
\]
where $Y$ is a semiresolution of $Z$~\cite{KoBa88}, and $\overline{Y}$ its normalization. 
Let $Q\in Z$ be a degenerate cusp. Then $\Delta = g^{-1}(Q)$ is either a nodal rational curve or a cycle of smooth rational curves. 
Since $\overline{C}$ is irreducible, then $C$ does not break in $\overline{Y}$ too. Now we consider all the possible cases.

\textbf{Case 1.} Suppose that $\Delta$ is a cycle of smooth rational curves. Since $C$ itself is smooth, $C$ must go through a 
vertex of the cycle, or intersect only one edge once. $\overline{Y} \la Y$ is generically $2-1$ on $\overline{C}$, since $Y-g^{-1}(Q)$ 
is semismooth~\cite{KoBa88}. And since $C$ does not break, $p^{-1} (\Delta)$ must remain a cycle. But then this cannot contract to a 
rational singularity. Therefore this case cannot happen.

\textbf{Case 2.} $\Delta$ is a nodal rational curve. Similarly as before, $\Delta$ remains a nodal rational curve in $\overline{Y}$, 
if $C$ does not break. 

Hence we proved that if $C$ does not break in $\overline{Z}$, $Z$ has no degenerate cusps and hence $\overline{Z}$ has exactly 
two singular points.

Now look at $\overline{Z} \stackrel{\overline{f}}{\la} T$. By~\cite[Proposition 4.26, 4.27]{KoBa88} it follows that the extended 
dual graph of $\overline{Z}$ must be \[
\circ \mbox{\noindent---} \cdots \mbox{\noindent---}\circ \mbox{\noindent---}\bullet \mbox{\noindent---}\circ 
\mbox{\noindent---} \cdots \mbox{\noindent---} \circ 
\]

Compute $\overline{C}^2$. Standard computations concerning cyclic quotient 
singularities show that \[
\overline{C}^2=-1 +\frac{\ts a}{\ts n}+\frac{\ts n-a}{\ts n}=0. \]
But this is impossible since $\overline{C}$ is contractible and hence it must be $\overline{C}^2 < 0$. 

Therefore we have shown that $C$ must break in the normalization of $Z$ to two smooth rational curves $C_1$ and $C_2$, 
and $\pi^{-1}(P)=\{P_1,\; P_2\}$, with $(P_1 \in \overline{Z}) \cong 1/n(1,a)$, and $(P_2 \in \overline{Z}) \cong 1/n(1,-a)$, 
where $P \in Z$ is the high index point. Moreover, since $C$ breaks, there is a degenerate cusp $Q \in Z$, and $\pi^{-1}(Q)=\{P\}=
C_1 \cap C_2$. We now want to understand the type of $P\in\overline{Z}$. 

In order to do this, we recall the following result by Shepherd-Barron.

\begin{lemma}[Lemma 1.3, ~\cite{Ba83}]\label{Barron}
Let $(P \in Z)$ be a degenerate cusp. Let $Y \stackrel{g}{\la}Z$ be the semi-resolution, and $\overline{Y} \stackrel{p}{\la} Y$, 
its normalization. Let $\Delta = g^{-1}(P)_{red}=\sum F_i$, the reduced exceptional locus of $g$, $E_i$ the strict transform of $F_i$ 
in $\overline{Y}$, $C$ the double curve of $Y$ and $B=p^{-1}(C)$.
Then
\begin{enumerate}
\item $F_i \cdot \Delta = E_i^2+2-E_i \cdot B $
\item $\text{mult}_P (Z)=\text{max}\{2,-\Delta^2\} $
\item $\text{embdim}_P (Z) = \text{max} \{3, -\Delta^2\}$.
\end{enumerate}
\end{lemma}

$Z$ corresponds to a divisorial neighborhood $Y \stackrel{f}{\la} X$. Since $Y$ is terminal and $P\in Z$ corresponds to an index one point of $Y$, 
it follows that $\text{embdim}_P(Z) \leq 4$. Hence from the previous lemma it follows that $-4 \leq \Delta^2 \leq -1$. 

By subadjunction we get that $p^{-1}(C)=\overline{C}_1+\overline{C}_2$. Now from part 1. of the previous lemma follows that \[
\Delta^2 = \sum_{i=1}^k E_i^2 + 2k-2. \label{Barron1}\]

We must now consider cases with respect to $\Delta^2$.

\textbf{Case 1.} $\Delta^2=-1$. I claim that in this case $Q \in \overline{Z}$ is a smooth point.

It is clear that there must be at least one $-1$ curve among the $E_i$'s. we want to find its position in $\Delta$. From~\ref{Barron1} 
follows that \[
\sum_{i=1}^k E_i^2=-2k+1. \]
 
From lemma~\ref{Barron} it follows that \[
F_i \cdot \Delta = E_i^2+2-(\overline{C}_1+\overline{C}_2)\cdot E_i= 1-(\overline{C}_1+\overline{C}_2)\cdot E_i  \]
If $E_i$ is not an edge of $\Delta$, then $F_i \cdot \Delta =1 \geq 0$, which is impossible. Hence the $-1$ curves must be edges of $\Delta$. 
Therefore, $\Delta$ is \[
\stackrel{-1}{\circ}\mbox{\noindent ---} \stackrel{-2}{\circ}\mbox{\noindent ---} \cdots \mbox{\noindent ---}\stackrel{-2}{\circ}\]
and hence $Q \in \overline{Z}$ is a smooth point. 

\textbf{Case 2.} $\Delta^2=-2$. In this case I claim that $Q \in \overline{Z}$ is an $A_m$ DuVal singularity. 

Again from lemma~\ref{Barron} follows that \[
\sum_{i=1}^k E_i^2=-2k. \]
If $E_i^2 \leq -2$ for all $i$, then $E_i^2=-2$ for all $i$, and hence $Q \in \overline{Z}$ is $A_k$. On the other hand if there is a 
$-1$ curve, then as before it must be an adge of $\Delta$. There are now two possibilities for $\Delta$. \[
\begin{array}{c}
\stackrel{-1}{\circ}\mbox{\noindent ---} \stackrel{-2}{\circ}\mbox{\noindent ---} \cdots \mbox{\noindent ---}\stackrel{-2}{\circ} 
\mbox{\noindent ---} \stackrel{-3}{\circ} \mbox{\noindent ---} \stackrel{-2}{\circ} 
\mbox{\noindent ---} \cdots \mbox{\noindent ---} \stackrel{-2}{\circ} \\
\stackrel{-1}{\circ}\mbox{\noindent ---} \stackrel{-2}{\circ}\mbox{\noindent ---} \cdots \mbox{\noindent ---}\stackrel{-2}{\circ} 
\mbox{\noindent ---} \stackrel{-3}{\circ} \mbox{\noindent ---} \stackrel{-2}{\circ} 
\mbox{\noindent ---} \cdots \mbox{\noindent ---} \stackrel{-2}{\circ} \mbox{\noindent ---} \stackrel{-3}{\circ} 
\mbox{\noindent ---} \stackrel{-2}{\circ} 
\mbox{\noindent ---} \cdots \mbox{\noindent ---} \stackrel{-2}{\circ} \mbox{\noindent ---} \stackrel{-1}{\circ}
\end{array}
\]
Both of them give an $A_m$ point as well. 

Working similarly we see that for the other two cases, the possible singularities for $Q \in \overline{Z} $ are:

\textbf{Case 3.} $\Delta^2=-3$, and hence $\text{mult}_P(Z)=3$. Then the dual graph in the minimal resolution of $Q \in \overline{Z}$ is \[
\stackrel{-2}{\circ}\mbox{\noindent ---} \cdots \mbox{\noindent ---}\stackrel{-2}{\circ} 
\mbox{\noindent ---} \stackrel{-3}{\circ} \mbox{\noindent ---} \stackrel{-2}{\circ}
\mbox{\noindent ---} \cdots \mbox{\noindent ---}\stackrel{-2}{\circ} \]

\textbf{Case 3.} $\Delta^2=-4$, and hence $\text{mult}_P(Z)=4$. Then the dual graph in the minimal resolution of $Q \in \overline{Z}$ is \[
\begin{array}{c}
\stackrel{-2}{\circ}\mbox{\noindent ---} \cdots \mbox{\noindent ---}\stackrel{-2}{\circ} 
\mbox{\noindent ---} \stackrel{-3}{\circ} \mbox{\noindent ---} \stackrel{-2}{\circ}
\mbox{\noindent ---} \cdots \mbox{\noindent ---}\stackrel{-2}{\circ}
\mbox{\noindent ---} \stackrel{-3}{\circ} \mbox{\noindent ---} \stackrel{-2}{\circ}
\mbox{\noindent ---} \cdots \mbox{\noindent ---}\stackrel{-2}{\circ} \\
\stackrel{-2}{\circ}\mbox{\noindent ---} \cdots \mbox{\noindent ---}\stackrel{-2}{\circ} 
\mbox{\noindent ---} \stackrel{-4}{\circ} \mbox{\noindent ---} \stackrel{-2}{\circ}
\mbox{\noindent ---} \cdots \mbox{\noindent ---}\stackrel{-2}{\circ} 
\end{array}
\]
In order to prove theorem~\ref{nonnormal1}.(b), it remains to describe the singularities with the above dual graph. 
\begin{lemma}
\begin{enumerate}
\item The singularity with dual graph \[
\stackrel{-2}{\circ}\mbox{\noindent ---} \cdots \mbox{\noindent ---}\stackrel{-2}{\circ} 
\mbox{\noindent ---} \stackrel{-3}{\circ} \mbox{\noindent ---} \stackrel{-2}{\circ}
\mbox{\noindent ---} \cdots \mbox{\noindent ---}\stackrel{-2}{\circ} \]
is the quotient singularity $1/m(1,-b)$, $(m,b)=1$, $b \not \equiv 1 \text{mod} m$, and $b | m+1$.
\item The singularity with dual graph \[
\stackrel{-2}{\circ}\mbox{\noindent ---} \cdots \mbox{\noindent ---}\stackrel{-2}{\circ} 
\mbox{\noindent ---} \stackrel{-4}{\circ} \mbox{\noindent ---} \stackrel{-2}{\circ}
\mbox{\noindent ---} \cdots \mbox{\noindent ---}\stackrel{-2}{\circ} \]
is the quotient singularity $1/m(1,-b)$, $(m,b)=1$, $b \not \equiv 1 \text{mod} m$, and $b | 2m+1$.
\item The singularity with dual graph \[
\stackrel{-2}{\circ}\mbox{\noindent ---} \cdots \mbox{\noindent ---}\stackrel{-2}{\circ} 
\mbox{\noindent ---} \stackrel{-3}{\circ} \mbox{\noindent ---} \stackrel{-2}{\circ}
\mbox{\noindent ---} \cdots \mbox{\noindent ---}\stackrel{-2}{\circ}
\mbox{\noindent ---} \stackrel{-3}{\circ} \mbox{\noindent ---} \stackrel{-2}{\circ}
\mbox{\noindent ---} \cdots \mbox{\noindent ---}\stackrel{-2}{\circ} \]
is the quotient singularity $1/m(1,-b)$, $(m,b)=1$, $b \not \equiv 1 \text{mod} m$, $1 \leq b < m$, $b | m+s+1$ for some $s \geq 1$, 
$b \geq 2s+1$ and $s+1 | b-s$.
\end{enumerate}
\end{lemma}
\begin{proof}
The proof of the lemma is just an exercise in continued fractions. I will only do $(2)$ here. The others are exactly similar with maybe 
more or less calculations.

So, suppose $P \in Z$ is a quotient singularity whose dual graph is the one described in $(2)$. Let $m, \; k$ be the number of 
$-2$'s to the left and to the  right of $-4$ respectively. Then 

\textbf{Claim.} \[
(P \in Z) \cong \frac{\ts 1}{\ts 2mk+3m+3k+4}(1,2mk+k+3m+1). \]

This can be checked by induction on $m$ and $k$ and the relation between the dual graph of a quotient singularity $1/n(1,a)$ 
and the continued fraction decomposition of $n/a$~\cite{Lam85}. 

Put $m=2mk+3m+3k+4$, and $b=2k+3$. Then $P \in Z \cong 1/m(1,-b)$. It is easy to see now that \[
2m+1=(2m+3)b \]
and therefore $b | 2m+1$, as claimed. Conversely, it is not difficult to see that any singularity $1/m(1,-b)$ such that $b | 2m+1$ 
has the dual graph that appears in $(2)$.  
\end{proof}

Next we want to find the conditions for $0 \in X$ to be an index 1 singularity, and $\Gamma$ smooth. To do this we will first show that 

\textbf{Claim:} $K_Y \cdot C =-1/n$. 

For this we will need the following easy lemma.

\begin{lemma}
Let $(P \in Z) \cong 1/n(1,a)$ be a quotient surface singularity, and let $C_1$, $C_2$ be two smooth proper curves intersecting each edge 
of the dual graph of $(P \in Z)$, i.e., the extended dual graph is \[
\stackrel{C_1}{\bullet}  \mbox{\noindent---} 
\stackrel{E_1}{\circ} \mbox{\noindent---} \cdots \mbox{\noindent---}\stackrel{E_k}{\circ} \mbox{\noindent---}\stackrel{C_2}{\bullet} \]
Then $C_1 \cdot C_2 = 1/n$.
\end{lemma}
\begin{proof}
Let $U \stackrel{f}{\la} Z$ be the minimal resolution of $(P\in Z)$, and $E_i$, $i=1, \ldots, k$ the exceptional curves. 
Suppose that $n/a=[b_1,\ldots,b_k]$, where $[b_1,\ldots,b_k]$ denotes the continued fraction decomposition of $n/a$. 
Then $E_i^2=-b_i$ and we can write \[
f^{\ast}C_1 = C_1^{\prime}+a_1E_1 + \cdots a_k E_k. \] 
Then $C_1 \cdot C_2 =a_k$. Intersecting with the $E_i$'s we get \[
\begin{array}{c}
0=1-b_1a_1+a_2 \\
0=a_1-b_2a_2+a_3 \\
\vdots \\
0=a_{k-1}-b_ka_k 
\end{array}
\]
Therefore we get that $a_{k-1}=b_ka_k$, $a_{k-2}=(b_kb_{k-1}-1)a_k$, $\dots$, $a_1= \beta a_k$, for some integer $\beta$. 
Hence from the above equations we see that $a_k = 1/ m$, for some integer $m$. But $a_1 = a/n $. Hence 
$a/n = \beta /m $, and therefore $am=\beta n$. Since $(a,n)=1$, it follows that $n | m$. But since $n a_1 \in \mathbb{Z}$, it follows 
from the above equations that $n a_i \in \mathbb{Z}$, for all $i$. Hence $na_k \in \mathbb{Z}$, and hence $m |n$. Therefore $m=n$, and 
the lemma follows.
\end{proof}

Return to the situation of theorems~\ref{nonnormal1}.~\ref{nonnormal2}. Assume that $\text{mult}_0\Gamma=\mu$ . 
Then as in the proof of lemma~\ref{multiplicity}, 
we see that $K_Y \cdot C = -\mu /\nu$, for some integer $\nu$. Now the index of $Y$ is $n$. Hence by adjunction it follows that \[
\omega_Z^{[n]} \otimes \sheaf_{C} = \sheaf_{\pone}(-n \mu /\nu).\]
By subadjunction it also follows that \[
\omega_{\overline{Z}}(C_1+C_2)= \pi^{\ast}\omega_Z \otimes \sheaf_{\overline{Z}}. \]
Hence \[
\omega_{\overline{Z}}^{[n]}(nC_1+nC_2)= \pi^{\ast}\omega_Z^{[n]} \otimes \sheaf_{\overline{Z}}. \]
Moreover, $\pi^{\ast}\omega_Z^{[n]} \otimes \sheaf_{C_i} \cong \omega_Z^{[n]} \otimes \sheaf_C \cong \sheaf_{\pone}(-n\mu /\nu)$, 
since $C_i \cong C$. Hence \[
\pi^{\ast}\omega_Z \cdot C_i = -\mu /\nu,\]
and hence by subadjunction again\[
C_i \cdot [ K_{\overline{Z}}+C_1+C_2] = -\mu /\nu.\]

We will now make a direct computation of the left hand side of the above equation. 
In the minimal resolution of $Z$, one of the curves $C_1$, $C_2$ must be a $-1$ 
curve. Suppose it is $C_1$. Then from standard computations in the theory of cyclic quotient singularities follows that \[
K_{\overline{Z}} \cdot C_1 = -1 -\frac{\ts a-n+1}{\ts n}-\frac{\ts m-b-m+1}{\ts m}=-\frac{\ts a+1}{\ts n}-\frac{\ts 1-b}{\ts m} \]
ans \[
C_1^2 = -1 +\frac{\ts a}{\ts n}+\frac{\ts m-b}{\ts m}=\frac{\ts a}{\ts n} - \frac{\ts b}{\ts m}.\]
Therefore, \[
C_i \cdot [ K_{\overline{Z}}+C_1+C_2] = -\frac{\ts 1}{\ts n} \]
and hence $\nu =n \mu $. Hence if $\Gamma$ is smooth, then $\mu =1$ and the claim as well as theorem~\ref{nonnormal2}.2 follows.

We now want to compute the multiplicity $\mu$ of $\Gamma$. 

\begin{lemma}\label{multiplicity1}
Let $C \subset Y \stackrel{f}{\la} X \ni 0 $ be a non-normal divisorial extremal neighborhood. Let $Z \in |\sheaf_Y|$ be the general member,
 $\overline{Z} \stackrel{\pi}{\la} Z$ be its normalization, and $T=f(Z)$. Write $K_{\overline{Z}}=\overline{f}^{\ast}K_T+x_1C_1+x_2C_2$. 
Let $E$ be the $f$-exceptional divisor and $\Gamma=f(E)$. Then the multiplicity $\mu$ of $\Gamma$ is\[
\mu = \frac{\ts x_1+x_2+2}{\ts n}.\]
\end{lemma}
\begin{proof}
Look at how $C$ degenerates in $Y$ (considered as a deformation of $Z$). As in lemma~\ref{multiplicity}, we see that $\mu$ copies of 
disjoint \pone's degenerate to $\nu C$, and we have already seen that $\mu = \nu / n$. Hence we need to find $\nu$. 
But $\nu$ is such that \[
[E \cap Z] = \nu C \]
as cycles. Moreover by adjunction we have that \[
\omega_Z=f^{\ast}\omega_T \otimes \sheaf_Y(E)\otimes \sheaf_Z . \]
Let $\mathcal{L}=\sheaf_Y(E)\otimes \sheaf_Z $. Then $\mathcal{L}^{-1}=I_{E\cap Z,Z}$. We want to find what is $\pi^{-1}(E \cap Z)$. 
By subadjunction we get that \[
\omega_{\overline{Z}}(C_1+C_2)=\pi^{\ast} \omega_Z = \overline{f}^{\ast}\omega_T \otimes \pi^{\ast} \mathcal{L},\]
and hence \[
\overline{f}^{\ast}\omega_T((x_1+1)C_1+(x_2+1)C_2)=\overline{f}^{\ast}\omega_T \otimes \pi^{\ast} \mathcal{L}. \]
Therefore, \[
\pi^{\ast} \mathcal{L} = \sheaf_{\overline{Z}}((x_1+1)C_1+(x_2+1)C_2), \]
and \[
[\pi^{-1}(E\cap Z)]=(x_1+1)C_1+(x_2+1)C_2. \]
Now consider the map \[
(x_1+1)C_1+(x_2+1)C_2 \la \nu C \]
and apply $\pi_{\ast}$ as cycle map. Then \[
\pi_{\ast} ((x_1+1)C_1+(x_2+1)C_2)=\nu C .\]
But $\pi_{\ast}[C_1]=\pi_{\ast}[C_2]=[C]$. Therefore $(x_1+x_2+2)[C]=\nu [C]$, and hence \[
\mu = \frac{\ts \nu}{\ts n}=\frac{\ts x_1+x_2+2}{\ts n}\]
and the lemma follows. 
\end{proof}

We now want to find the condition for $0 \in T$ to be DuVal and hence $0 \in X$ cDV. It is clear that $0 \in T$ is DuVal 
if and only if 
\begin{enumerate}
\item $x_1, \; x_2 \in \mathbb{Z}$, and
\item $K_{\overline{Z}}-x_1C_1-x_2C_2 \; \in \text{Pic}(\overline{Z})$
\end{enumerate}
We will now interpret these conditions. Condition $(1)$ above 
means that there are integer solutions of the system of equations \[
\begin{array}{c}
x_1C_1^2+x_2C_1\cdot C_2 = K_{\overline{Z}} \cdot C_1 \\
x_1 C_1 \cdot C_2 +x_2 C_2^2=K_{\overline{Z}} \cdot C_2
\end{array} \]
By lemma~\ref{multiplicity1}, $\Gamma$ is smooth iff $x_1+x_2=n-2$. This together with the above system of equations give that 
\[
x_1=\frac{\ts K_{\overline{Z}} \cdot C_1-(n-2)C_1 \cdot C_2}{\ts C_1^2-C_1\cdot C_2}=
\frac{\ts K_{\overline{Z}} \cdot C_2-(n-2)C_2^2 }{\ts C_1 \cdot C_2 - C_2^2} \in \mathbb{Z}\]
We also know that $[K_{\overline{Z}}+C_1+C_2]\cdot C_i = -1/n$. $i=1,\; 2$. Therefore, \[
K_{\overline{Z}} \cdot C_2-(n-2)C_2^2 = -1/n-C_1 \cdot C_2 - (n-1)C_1^2. \]
Hence \[
\frac{\ts K_{\overline{Z}} \cdot C_1-(n-2)C_1 \cdot C_2}{\ts C_1^2-C_1\cdot C_2}=
-\frac{\ts 1/n + C_1 \cdot C_2 +(n-1)C_2^2}{\ts C_1 \cdot C_2 - C_2^2} 
\]
subtracting $n-1$ from both sides we get \[
\frac{\ts K_{\overline{Z}} \cdot C_1+C_1\cdot C_2 -(n-1)C_1^2}{\ts C_1^2-C_1\cdot C_2}=
-\frac{\ts 1/n + n C_1 \cdot C_2 }{\ts C_1 \cdot C_2 - C_2^2} 
\]
Use that $[K_{\overline{Z}}+C_1+C_2]\cdot C_i = -1/n$ again to get that \[
\frac{\ts 1+n^2C_1^2}{\ts C_1^2-C_1\cdot C_2}=
\frac{\ts 1+n^2 C_1 \cdot C_2 }{\ts C_1 \cdot C_2 - C_2^2} 
\]
This in return gives that \[
(C_1-C_2)^2=n^2[(C_1 \cdot C_2)^2 -C_1^2 C_2^2], \]
and theorem~\ref{nonnormal1}.3 follows.

We also need that $x_1 \in \mathbb{Z}$. As before we see that $K_{\overline{Z}}\cdot C_1 = -a/n-1/n+b/m-1/m$, $C_1^2=a/n-b/m$, 
and $C_1 \cdot C_2=1/m$. Then a simple calculation shows that 
\begin{equation}
x_1=\frac{\ts K_{\overline{Z}} \cdot C_1-(n-2)C_1 \cdot C_2}{\ts C_1^2-C_1\cdot C_2}=-1 - \frac{\ts m+n^2}{\ts am-(b+1)n}.
\end{equation}
Therefore $a,\;  b, \;  m, \;  n$, must satisfy 
\begin{equation}
\frac{\ts m+n^2}{\ts am-(b+1)n} \in \mathbb{Z}. 
\end{equation}
We now want to find the conditions so that $K_{\overline{Z}}-x_1C_1-x_2C_2 \; \in \text{Pic}(\overline{Z})$. At first we notice that 
this condition need only be satisfied at $P_1$, $P_2$. The reason is that $Q$ is over a degenerate cusp which corresponds to an 
index terminal 1 singularity of $Y$ and hence this condition is always satisfied there. 

The following result is useful.

\begin{lemma}
Let $(0 \in Z)\cong 1/n(1,a) $ be a cyclic quotient surface singularity, and $C$ a smooth curve that intersects an edge of the 
dual graph of $0 \in Z$. Then $K_Z - m C \in \text{Pic}(Z)$, where $m >0$ is a positive integer such that $n | (am+a+1)$.
\end{lemma}
\begin{proof}
Let $U \stackrel{f}{\la} Z$ be the minimal resolution and $E_i$ the exceptional divisors. Then we can write \[
\begin{array}{c}
K_U=f^{\ast}K_Z+ \sum_{i=1}^k \beta_i E_i \\
f^{\ast}C=C^{\prime}+\sum_{i=1}^k \gamma_i E_i
\end{array} \]
We must find an integer $m>0$ such that $\beta_i+m\gamma_i \in \mathbb{Z}$. 
But we also know that $\beta_1=(a-n+1)/n$, and $\gamma_1=a/n$. Now as in the proof of lemma 2.11 we see that it is sufficient to find an $m$ 
such that $\beta_1+m \gamma_1 \in \mathbb{Z}$. But this is equivalent to the condition of the lemma.
\end{proof}
Now if $0 \in T$ is DuVal, then $K_{\overline{Z}}-x_1C_1-x_2C_2$ is Cartier. We have already said that this condition is in any case 
satisfied at $Q \in \overline{Z}$ since it lies over a degenerate cusp that corresponds to an isolated index 1 point of $Y$. So 
we only need to find conditions at $P_1,\; P_2$. By lemma~\ref{multiplicity1}, and since $\Gamma$ is smooth, $x_1+x_2=n-2$. Therefore we need 
to check when $K_{\overline{Z}}-x_1C_1-(n-2-x_1)C_2$ is Cartier. Check at $P_1$ first. By the previous lemma this happens iff 
$n|ax_1+a+1$. Similarly, at $P_2$ this happens iff $n |(n-a)(n-2-x_1)+n-a+1$, which is equivalent to $n|ax_1+a+1$. 
But by $(2)$ follows that \[
ax_1+a+1=- n\frac{\ts b+1+an}{\ts am-(b+1)n}. \]
Hence $n | ax_1+a+1$ iff \[
\frac{\ts b+1+an}{\ts am-(b+1)n} \in \mathbb{Z}. \] 
This concludes the proof of theorem~\ref{nonnormal1}. 

Theorem~\ref{nonnormal2}.1, .4 are proved exactly the same way as theorem~\ref{normal2}. It only remains to show $(3)$. 

Pick a line $0 \in l_1 \subset T$ such that it intersects the edge of the dual graph of $0 \in T$, and that 
$\overline{l}_1=\overline{f}_{\ast}^{-1}l$ goes through $P_1$. Then if $0 \in T$ is $A_{d-1}$, $dl_1$ is Cartier and hence the index 
of $\overline{f}^{\ast}l_1$ is $d$. 

Suppose that \[
\overline{f}^{\ast}l_1=\overline{l}_1+a_1C_1+a_2C_2 \]
Intersect with $C_1, \; C_2$. Then \[
\begin{array}{c}
a_1C_1^2+a_2C_1 \cdot C_2 = -C_1 \cdot \overline{l}_1=-1/n \\
a_1 C_1 \cdot C_2 +a_2C_2^2=0 
\end{array} \]
Therefore we find that \[
\begin{array}{cc}
a_1=n\frac{\ts C_2^2}{\ts (C_1-C_2)^2}, & a_2 = -n\frac{\ts C_1 \cdot C_2}{\ts (C_1-C_2)^2}=-n\frac{\ts 1}{\ts m(C_1-C_2)^2} 
\end{array} \]
Look at $P_2$. $\text{index}_{P_2}(\overline{f}^{\ast}l_1)=\text{index}_{P_2}(n/m(C_1-C_2)^2)C_2$. $nC_2$ is Cartier at $P_2$. 
Hence if $(C_1-C_2)^2=\nu / k$, then $\text{index}_{P_2}(\overline{f}^{\ast}l_1)=m \nu / (k,m \nu)$. 

We must now compute the index of $C_1-C_2$. At $P_1$, $P_2$ it is $n$. Check now at $Q$. The extended dual graph of $Q \in \overline{Z}$ is 
\[
\bullet \mbox{\noindent ---} \circ\mbox{\noindent ---}\cdots \mbox{\noindent ---} \circ\mbox{\noindent ---}\bullet \]
Let $U \stackrel{g}{\la} Q \in \overline{Z}$ be the minimal resolution, and let $E_i$ be the exceptional curves. Then \[
\begin{array}{c}
g^{\ast}C_1=C_1^{\prime}+\frac{\ts m-b}{\ts m}E_1 + \cdots \\
g^{\ast}C_2=C_2^{\prime}+\frac{\ts 1}{\ts m}E_1 + \cdots
\end{array} \]
Therefore \[
g^{\ast}(C_1-C_2)=C_1^{\prime}-C_2^{\prime}+\frac{\ts m-b-1}{\ts m}E_1 + \cdots 
\]
and hence \[
\text{index}_Q (C_1-C_2)=\frac{\ts m}{\ts (m,b+1)} \]
Therefore \[
k=\text{index}(C_1-C_2)=lcm(n,\frac{\ts m}{\ts (m,b+1)} )\]
I now claim that in fact $k=mn/(m,b+1)$.This will follow if we prove that $(n,m/(m,b+1))=1$. Suppose that a prime $p$ divides 
$n$ and $m/(m,b+1)$. We have already shown that $am-(b+1)n$ divides $b+1+an$. Hence $p |(b+1)$ and of course $p |m$. Hence $p=1$.  

Now let $k(C_1-C_2)^2=\nu \in \mathbb{Z}$. Then from the previous discussion it follows that \[
\text{index}_{P_2}(\overline{f}^{\ast}l_1)=\frac{\ts m \nu}{\ts (k,m \nu)}\]

We now want to check at $P_2$. To do this pick another line $0 \in l_2 \subset T$ such that it intersects one end of the dual graph 
of $0 \in T$ and that $\overline{l}_2$ goes through $P_2$. Then $l_1+l_2=K_T$ and it is Cartier. 

\textbf{Claim:} \[
\overline{f}^{\ast}(l_1+l_2)=\overline{l}_1+\overline{l}_2+(x_1+1)C_1+(x_2+1)C_2. \]
Then at $P_1$ it follows that $\overline{l}_1+(x_1+1)C_1$ is Cartier, and hence we need to find the index at $P_1$ of 
$(a_1-x_1-1)C_1$. But a simple calculation shows that $a_1+a_2=x_1+1$. Hence $(a_1-x_1-1)C_1=-a_2C_1=-n/m(C_1-C_2)^2 C_1$. 
As before we conclude again that \[
\text{index}_{P_1}(\overline{f}^{\ast}l_1)=\frac{\ts m \nu}{\ts (k,m \nu)}\]
The same method applies to $Q$ also and we get theorem~\ref{nonnormal2}.3. 

Now to prove the claim. Suppose that \[
\overline{f}^{\ast}(l_1+l_2)=\overline{l}_1+\overline{l}_2+b_1C_1+b_2C_2 \]
and let $Y \stackrel{f}{\la} X$ the corresponding threefold contraction. 
Let $S$ be the general section of $X$ through $\Gamma$. Then $S \cap T =l_1+l_2$. 
Moreover by assumption, $S_Y \cong S$. Therefore $S_Y \cap Z = l_1^{\prime}+l_2^{\prime}$. 
But $S_Y \in |-K_Y|$ and therefore $S_Y \cdot Z =-K_Y \cdot Z = - E \cdot Z$. 
Now by using subadjunction and adjunction as before we conclude that \[
\overline{f}^{\ast}(l_1+l_2)=\overline{l}_1+\overline{l}_2+(x_1+1)C_1+(x_2+1)C_2 \]
as claimed.

\end{proof}

The next example shows that it is possible that $Q\in \overline{Z}$ is singular and in the resulting contraction $Y \la X$, 
$\Gamma$ is smooth.

\textbf{Example:}  Consider the configuration \[
\stackrel{-2}{\circ} \mbox{\noindent ---} \stackrel{-2}{\circ} \mbox{\noindent ---}
\stackrel{-2}{\circ} \mbox{\noindent ---} \stackrel{-3}{\circ} \mbox{\noindent ---} 
\stackrel{-2}{\bullet} \mbox{\noindent ---} \stackrel{-3}{\circ} \mbox{\noindent ---}
\stackrel{-1}{\bullet} \mbox{\noindent ---} 
\stackrel{-2}{\circ} \mbox{\noindent ---} \stackrel{-5}{\circ} 
\]
Contract all the curves except the ones marked by a solid circle. Then we get a map $\overline{Z} \stackrel{f}{\la} T$, such that $0 \in T$ is 
an $A_5$ singularity, and $\overline{Z}$ has exactly three singular points $P_1 \in C_1$, $P_2\in C_2$ and $Q \in C_1 \cap C_2$. 
Moreover it is easy to see that $(P_1 \in \overline{Z}) \cong 1/9 (1,5)$, $(P_2 \in \overline{Z}) \cong 1/9(1,-5)$, 
and $Q \in \overline{Z}) \cong 1/3(1,1)$. 

Now it is not difficult to see that $C_1^2=-1+5/9+1/3=-1/9$, $C_2^2=-2+4/9+1/3=-11/9$, $C_1 \cdot C_2 =1/3$, 
$K_{ \overline{Z}}\cdot C_1=-1/3$ and $K_{ \overline{Z}}\cdot C_2=7/9$. 
Then in the notation of the proof of theorem~\ref{nonnormal1}, $x_1=6$ and $x_2=1$. Therefore \[
\text{mult}_0 (\Gamma)=\frac{\ts 6+1+2}{\ts 9}=1 \]
and hence $\Gamma$ is smooth.

%% file: section4.tex
\section{Semistable neighborhoods}

In this section we will classify the semistable divisorial extremal neighborhoods $C \subset Y \stackrel{f}{\la} X \ni 0$. 
Let $E$ be the $f$-exceptional divisor, $\Gamma =f(E)$, and $S$ the general section of $X$ through $\Gamma$. 
Recall that such a neighborhood is called semistable if $C \subset S_Y = f_{\ast}^{-1}S$. 

The following example shows that such neighborhoods do exist.

\textbf{Example:}  Consider the configuration \[
\stackrel{-3}{\circ} \mbox{\noindent ---}
\stackrel{-2}{\circ} \mbox{\noindent ---} \stackrel{-3}{\circ} \mbox{\noindent ---} 
\stackrel{-1}{\bullet} \mbox{\noindent ---} \stackrel{-2}{\circ} \mbox{\noindent ---}
\stackrel{-5}{\circ}  
\]
Contract all curves except the $-1$ curve to obtain a morphism $C \subset Z \stackrel{f}{\la} T \ni 0$, where $0 \in T$ is an $A_2$ 
DuVal singularity, and $Z$ has exactly two singular points $P, \; Q$. Moreover by the construction it is clear that 
$(P \in Z) \cong 1/12(1,5)$, and $(Q \in Z) \cong 1/9(1,5)$. $Z$ has index $2$ at $P$ and $3$ at $Q$. These singularities are $T$-singularities, 
and therefore they admit $\mathbb{Q}$-Gorenstein terminal smoothings. Take one, say $Y$. Then $f$ extends to the family and we 
get a semistable divisorial extremal neighborhood $Y \stackrel{f}{\la} X$. Let $E$ be the exceptional divisor and $\Gamma=f(E)$.
As in previous examples, it is easy to compute that 
$K_Z \cdot C=-1/6$, and $C^2=-1/36$. Then by lemma~\ref{multiplicity}, 
$mult(\Gamma)=-(K_Y \cdot C)^2 / C^2 = 1$, and hence $\Gamma$ is a smooth curve.

\textbf{Remark:} At this point I would like to point out that semistable extremal neighborhoods cannot be constructed from 
the diagram of theorem~\ref{normal2}.1. The reason is that if this was possible then $S_Y \cong S$, which is not the case here. 
However, they can be constructed from the diagram~\cite{Tzi02}
\[
\xymatrix{ 
V \ar[d]_{\nu} \ar@{-->}[rr]^{\phi} & & V^{\prime} \ar[d]^{\pi}\\
W \ar[dr]_g &   & Y \ar[dl]^f \\
            & X & }
\]
but flips will always appear in $\phi$.

Let $0 \in T$ be the general section of $X$ through $0$. As before we want to understand the singularities of $Z=f^{\ast}T$, 
and then treat the problem of classification of semistable extremal neighborhoods as a problem of deformations of certain surface singularities. 

\begin{lemma}\label{sms-singularities}
Let $C \subset Y \stackrel{f}{\la} X \ni 0$ be a semistable divisorial extremal neighborhood, and let $0 \in T$ be the 
general section of $X$ through $0$. Then $Z=f^{\ast}T$ is normal and its singular locus is exactly two T-singularities of index 
bigger than one.
\end{lemma}
\begin{proof}
Let $S$ be the general section of $X$ through $\Gamma$. By the proof of lemma~\ref{singularities1}, the pair $(Y,S_Y+Z)$ is 
log-canonical. Since the neighborhood is semistable, we have that $C \subset S_Z$. Therefore the pair $(Y,Z)$ is $klt$. Hence 
$Z$ is normal~\cite[Proposition 5.51]{Ko-Mo98} and the high index points are by lemma~\ref{singularities1} $T$-singularities. 
Moreover, by~\cite[Theorem 2.2.4]{Ko-Mo92}, $Y$ has exactly two singular points and they both have index bigger than one. 
Hence also $Z$ has exactly two high index points. It remains to show that $Z$ has no other singular points. 

Again from~\cite[Theorem 2.2.4]{Ko-Mo92} follows that the extended dual graph of $S_Y$ is 
\[
\circ \mbox{\noindent---} \cdots \mbox{\noindent---}\circ \mbox{\noindent---}\stackrel{C}{\bullet} \mbox{\noindent---}\circ 
\mbox{\noindent---} \cdots \mbox{\noindent---} \circ 
\]
$0\in S$ can be written as $xy-z^{n+1}=0$. Let $l_1$, $l_2$ be the lines $x=z=0$ and $y=z=0$ respectively. 
Now consider the cycle $\Delta=l_1^{\prime}+l_2^{\prime}+C$. It is not difficult to see by examining the dual graph, 
that $C \cdot \Delta =0$, and in fact $\Delta =f^{\ast}(l_1+l_2)$, and in particular it is Cartier. 

\textbf{Claim:} There is a surface $H \subset Y$ such that $H \cap S_Y = \Delta$, and $C \cdot H =0$.

Then $S_Y$ is Cartier at any index one point of $Y$, and hence $\Delta = S_Y \cap H$ is also Cartier in $H$. But since $\Delta$
 is smooth  away from the singular points of $Y$, it follows that $H$ must be smooth as well. Then $T=f(H)$ is a section of $X$ 
through $0$ and we can take $Z$ to be $H$. 

Now to prove the claim. $\Delta \in |\sheaf_{S_Y}|$. Since $S_Y \in |-K_Y|$, there is an exact sequence \[
0 \la \sheaf_{Y}(K_Y) \la \sheaf_Y \la \sheaf_{S_Y} \la 0 \]
By Kawamata-Viehweg vanishing it follows that the following sequence is exact \[
0 \la H^0(\sheaf_Y(K_Y)) \la H^0(\sheaf_Y) \la H^0(\sheaf_{S_Y})\la H^1(\sheaf_Y(K_Y))=0 \]
Therefore there exist an element $H \in |\sheaf_Y|$ with the required properties. Now the claim and hence the lemma too follows.
\end{proof}

We may now use the deformation theory of $T$-singularities to classify semistable neighborhoods.

\begin{theorem}\label{semistable1}
Semistable divisorial extremal neighborhoods with smooth center are in one to one correspondence with $\mathbb{Q}$-Gorenstein 
smoothings of germs of surfaces $C \subset Z$, with $C \cong \pone$, with the following properties
\begin{enumerate}
\item $Z$ has exactly two singular points $P$ and $P^{\prime}$. Moreover, \[
\begin{array}{cc}
(P \in Z)\cong \frac{\ts 1}{\ts n^2 d}(1,and-1), & 
(P^{\prime} \in Z)\cong \frac{\ts 1}{\ts {n^{\prime}}^2 d^{\prime}}(1,a^{\prime}n^{\prime}d^{\prime}-1)
\end{array} \]
\item $(n,n^{\prime})=(a,a^{\prime})=1$.
\item $d | n^{\prime}$, and $d^{\prime} | n$.
\item $nn^{\prime}-an^{\prime}-a^{\prime}n=-1$.
\item $dd^{\prime}={n^{\prime}}^2d^{\prime}+n^2d-nn^{\prime}dd^{\prime}$.
\end{enumerate}
\end{theorem}

\begin{proof}
Let $C \subset Y \stackrel{f}{\la} X \ni 0$ be a semistable divisorial extremal neighborhood. Let $0 \in T$ be the general 
section of $X$ through $0$ and $Z=f^{\ast}T$.

Part $(1)$ follows immediately from corollary~\ref{singularities2} and lemma~\ref{sms-singularities}. 

To see $(2)$. By~\cite[Corollary 1.10]{Mo88}, it follows that $(n,n^{\prime})=1$ iff $\text{Cl}^{sc}(Y)$ is torsion free. 
But we have already seen in the proof of theorem~\ref{normal1} that $X$ has index one iff this happens. Therefore $(n,n^{\prime})=1$. 

To continue we now compute $K_Z \cdot C$ and $C^2$. Standard computations in cyclic quotient singularities give that 
\begin{eqnarray}
K_Z \cdot C=1-\frac{a}{n}-\frac{a^{\prime}}{n^{\prime}}=\frac{nn^{\prime}-an^{\prime}-a^{\prime}n}{nn^{\prime}} \\
C^2=\frac{nn^{\prime}dd^{\prime}-{n^{\prime}}^2d^{\prime}-n^2d^{\prime}}{n^2{n^{\prime}}^2dd^{\prime}}
\end{eqnarray}
Now arguing as in the proof of lemma~\ref{multiplicity} we see that $K_Y \cdot C = -1/\nu$, for some $\nu$. Therefore 
by adjunction it follows that $K_Z \cdot C= K_Y \cdot C = -1/ \nu$, and therefore \[
\frac{1}{K_Z \cdot C} = \frac{nn^{\prime}}{nn^{\prime}-an^{\prime}-a^{\prime}n} \in \mathbb{Z} . \]
Let $p$ be a prime dividing $nn^{\prime}-an^{\prime}-a^{\prime}n$. Then $p|n$ or $p| n^{\prime}$. But $(a,a^{\prime})=(n,n^{\prime})=1$. 
If $p | n$ then $p | an^{\prime}$, and hence $p=1$. Similarly if $p|n^{\prime}$. Therefore we conclude that 
$nn^{\prime}-an^{\prime}-a^{\prime}n=-1$, and part $(4)$ follows. 

From lemma~\ref{multiplicity} it follows that the center $\Gamma$ is smooth iff \[
\frac{(K_Z \cdot C)^2}{C^2}=-1 \]
and hence iff 
\begin{equation}
dd^{\prime}={n^{\prime}}^2d^{\prime}+n^2d-nn^{\prime}dd^{\prime}, 
\end{equation}
and part $(5)$ follows.

It remains to show $(3)$. It is clear that $\text{index}(Z)=\text{lcm}(n,n^{\prime})=nn^{\prime}$. Moreover, 
$\text{index}_P(C)=n^2d$, and $\text{index}_{P^{\prime}}(C)={n^{\prime}}^2d^{\prime}$. It is also clear from $(4), \; (5),\; (6)$, 
that \[
K_Z=f^{\ast}K_T+nn^{\prime}C.\]
Hence at $P$, $n^2n^{\prime}C$ is Cartier. Therefore $n^2d|n^2n^{\prime}$, and hence $d|n^{\prime}$. Similarly checking at $P^{\prime}$ 
we get that $d^{\prime} | n$. 
\end{proof}

Now let $C \subset Z$ be a germ as in the theorem, and let $Y \stackrel{f}{\la} X$ be the corresponding extremal neighborhood. 
The next theorem gives information about the singularities of $X$.
\begin{theorem}\label{semistable2}
Let $C \subset Z$ be a surface germ as in the previous theorem, and let $Y \la X$ be the corresponding extremal neighborhood. 
Then
\begin{enumerate}
\item The general section $T$ of $X$ through $0$ is an $A_{dd^{\prime}-1}$ DuVal singularity.
\item The general section $S$ of $X$ through the center $\Gamma$ is an $A_{nk+n^{\prime}k^{\prime}-1}$ DuVal singularity, 
where $k,\; k^{\prime}$ are the axial multiplicities of $P\in Y$, $P^{\prime}\in Y$.
\item Let $E_i$ be the exceptional curves in the minimal resolution of $S$. Then $\Gamma$ intersects $E_s$, with 
$s=nk^{\prime}/d^{\prime}-n^{\prime}k/d+kn$.
\end{enumerate}
\end{theorem}

\begin{proof}
To show $(2)$ we will follow the same method as in the proof of theorem~\ref{nonnormal2}. 

$0 \in T$ is of the form $xy-z^{m+1}=0$. Now let $l_1$, $l_2$ be the lines $x=z=0$, and $y=z=0$, respectively. 
Then they both have index $m+1$ in $T$. Let $l_1^{\prime}$, $l_2^{\prime}$ be their birational transforms in $Z$. 
Then $P \in l_1^{\prime}$, and $P^{\prime} \in l_2^{\prime}$. Then the index of $f^{\ast}l_i$ is $m+1$ and this is what we want 
to find. Suppose that $f^{\ast}l_1=l_1^{\prime}+\delta C$. From lemma 4.11 it follows that $l_1^{\prime}\cdot C=1/n^2d$, and 
$l_2^{\prime}\cdot C = 1/{n^{\prime}}^2d^{\prime}$. Hence intersecting with $C$ we find that 
\begin{eqnarray*}
f^{\ast}l_1=l_1^{\prime}+\frac{{n^{\prime}}^2}{d}C \\
f^{\ast}l_2=l_2^{\prime}+\frac{n^2}{d^{\prime}}C 
\end{eqnarray*}
$l_1+l_2$ is Cartier in $T$ and hence \[
f^{\ast}(l_1+l_2)=l_1^{\prime}+l_2^{\prime}+(\frac{{n^{\prime}}^2}{d}+\frac{n^2}{d^{\prime}})C \]
is also Cartier. Checking at $P$ we see that $(l_1^{\prime}+{n^{\prime}}^2/d C)+n^2/d^{\prime}C$ is also Cartier. 
Therefore the index of $f^{\ast}(l_1)$ is $dd^{\prime}$ and part $(1)$ follows. 

Part $(2)$ follows from~\cite[Theorem 2.2.4]{Ko-Mo92}. 

Now let $E$ be the $f$-exceptional divisor and $\Gamma=f(E)$ it's center that by assumption is a smooth curve. 
We want to find the position of $\Gamma$ in the fundamental cycle of the general section $S$ of $X$ through $\Gamma$. 
For this we need the following easy result.

\begin{lemma}
Let $(0 \in C \subset T)$ be the germ of an $A_n$ DuVal singularity alond a proper curve $C$. Let $U \stackrel{f}{\la} T$ be 
the minimal resolution and $E_1, \ldots, E_n$ the exceptional curves. Let $\Gamma$ be a smooth curve through $0$ such that 
$f_{\ast}^{-1} \Gamma$ intersects $E_k$. Then\[
C \cdot \Gamma =\frac{k}{n+1}
\]
\end{lemma}
\begin{proof}
It is easy to see that \[
f^{\ast}C=C^{\prime}+\sum_{i=1}^n \frac{i}{n+1}E_i \]
and therefore $C \cdot \Gamma =k/(n+1)$.
\end{proof}
Now by~\cite[Theorem 2.2.4]{Ko-Mo92}, the extended dual graph of $S_Y$ is 
\begin{equation}
\stackrel{E_1}{\circ} \mbox{\noindent ---}
\cdots  \mbox{\noindent ---} \stackrel{E_{nk-1}}{\circ} \mbox{\noindent ---} 
\stackrel{C}{\bullet} \mbox{\noindent ---} \stackrel{E_{nk+1}}{\circ} \mbox{\noindent ---} \cdots \mbox{\noindent ---}
\stackrel{E_{n^{\prime}k^{\prime}+nk-1}}{\circ} 
\end{equation}
Now we want to find the divisorial part of $S_Y \cap E$. Suppose that \[
(S_Y \cap E)_{div} = \Gamma + \delta C \]
for some $\delta$. Suppose that $\Gamma$ intersects $E_m$, with $m < nk$. 
 Pick a line $l$ through $0$ in $S_Y$ such that in the minimal resolution it intersects $E_{n^{\prime}k^{\prime}+nk-1}$. 

Now compute $l \cdot E$ in two ways.
\[
l \cdot E = (l \cdot (E|S_Y))_{S_Y} = (l \cdot (\Gamma +\delta C))_{S_Y}=\delta (l \cdot C)_{S_Y}= \frac{\delta}{n^{\prime}k^{\prime}} \]
and \[
l \cdot E=(l \cdot E|Z)_Z=(l \cdot (nn^{\prime}C))_Z=nn^{\prime}\frac{1}{{n^{\prime}}^2d^{\prime}}=\frac{n}{n^{\prime}d^{\prime}} \]
Here we have used that $E|Z=nn^{\prime}C$. This follows from adjunction in $Y$ since $K_Z =f^{\ast}K_T +nn^{\prime}C$, and lemma 4.11. 
Therefore $\delta = nk^{\prime}/d^{\prime}$. 

We now want to compute $(C^2)_{S_Y}$. Let $U \stackrel{g}{\la} S_Y$ be the minimal resolution. Then from $(7)$ it follows that \[
g^{\ast}C = C^{\prime}+\sum_{i=1}^{nk-1} \frac{i}{nk}E_i + \sum_{j=1}^{n^{\prime}k^{\prime}-1} \frac{n^{\prime}k^{\prime}-j}{n^{\prime}k^{\prime}}E_j \]
Therefore, \[
(C^2)_{S_Y}=-2+\frac{nk-1}{nk}+\frac{n^{\prime}k^{\prime}-1}{n^{\prime}k^{\prime}}=-\frac{1}{nk}-\frac{1}{n^{\prime}k^{\prime}} \]

Now 
\begin{gather*}
-\frac{1}{nn^{\prime}}=C \cdot K_Y=C \cdot E=(C \cdot E|{S_Y})_{S_Y}=(C \cdot (\Gamma + \frac{nk^{\prime}}{d^{\prime}}C))|_{S_Y}= \\
=(C \cdot \Gamma)_{S_Y} -\frac{nk^{\prime}}{d^{\prime}}\frac{nk+n^{\prime}k^{\prime}}{nn^{\prime}kk^{\prime}} 
\end{gather*}
Now by the previous lemma, \[
(C \cdot \Gamma)_{S_Y}=\frac{m}{nk}\]
and therefore 
\begin{gather*}
m=\frac{nk^{\prime}}{d^{\prime}}\frac{nk+n^{\prime}k^{\prime}}{n^{\prime}k^{\prime}}-\frac{k}{n^{\prime}}= 
\frac{n}{d^{\prime}}k^{\prime}+\frac{n^2-d^{\prime}}{n^{\prime}d^{\prime}}k
\end{gather*}
But by theorem~\ref{semistable1}.5, $dd^{\prime}=n^2d+{n^{\prime}}^2d^{\prime}-nn^{\prime}dd^{\prime}$, and hence $(n^2-d^{\prime})d=
-n^{\prime}d^{\prime}(n^{\prime}-nd)$. Therefore, \[
m=\frac{n}{d^{\prime}}k^{\prime}-\frac{n^{\prime}-nd}{d}k=\frac{n}{d^{\prime}}k^{\prime}-\frac{n^{\prime}}{d}k+kn \]
as claimed.
\end{proof}

We can now relate extremal neighborhoods with the solutions of certain Pell equations.
\begin{corollary}\label{Pell}
Fix an $A_{k-1}$ DuVal singularity $0 \in T$. Then semistable divisorial extremal neighborhoods $C \subset Y \stackrel{f}{\la} X \ni 0$, 
with smooth center, and such that the general section of $X$ through $0$ is $T$, are in one to one correspendence with the integer solutions $(x,y)$, 
with $y \neq 0$, of the Pell equations \[
\frac{k}{d}x^2-d(k-4)y^2=4 \]
 for all $d | k$. In particular, if $k=p$ is prime, then such semistable neighborhoods are in one to one correspondence with the integer 
solutions of the equation $x^2-p(p-4)y^2=4$.
\end{corollary}
\textbf{Remark:} In general it is not known that the above Pell equations always have solutions. However, if $p(p-4)$ is square free and positive, 
then there are infinitely many and hence infinitely many corresponding semistable neighborhoods. But it is possible that for special 
values of $p$, there may be no solutions and hence no corresponding semistable neighborhoods. Unfortunately I do not know of any such examples. 
The case $p=2, \; 3$ are particularly simple and we will treat them later. 
\begin{proof}
From theorem~\ref{semistable1}.5, $dd^{\prime}=n^2d+{n^{\prime}}^2d^{\prime}-nn^{\prime}dd^{\prime}$. Moreover, $d|n^{\prime}$, and $d^{\prime}|n$. 
Hence, $n^{\prime}=m^{\prime}d$, and $n=md^{\prime}$. Hence 
\begin{equation}
1={m^{\prime}}^2d+m^2d^{\prime}-mm^{\prime}dd^{\prime}. 
\end{equation}
By theorem~\ref{semistable2}.1, we know that $k=dd^{\prime}$. Hence putting $d^{\prime}=k/d$ and completing the square in $(8)$, 
we find that \[
\frac{k}{d}(2m-dm^{\prime})^2-d(k-4){m^{\prime}}^2=4 \]
and therefore we get the correspondence between semistable extremal neighborhoods and integer solutions of the Pell equation 
$k/d x^2-d(k-4)y^2=4$ as claimed. Moreover, the correspondence is given by $x=2n/d^{\prime}-n^{\prime}$, and $y=n^{\prime}/d$. 
Clearly now $y$ cannot 
be zero.
\end{proof}
\begin{proof}[Proof of Theorem~\ref{deformation}]
Fix data as in theorem~\ref{deformation}. Let $C \subset Y \stackrel{f}{\la} X \ni 0$ be a divisorial extremal neighborhood. 
Then let $Z=f^{\ast}T$. By~\cite[Proposition 11.4]{Ko-Mo92}, there is a map 
\begin{equation}
Def(Z) \la Def(T).
\end{equation} 
By theorems~\ref{normal1},~\ref{normal2},
~\ref{nonnormal1},~\ref{nonnormal2},~\ref{semistable1},~\ref{semistable2}, there are finitely many only surface germs $C \subset Z$, 
so that the assumptions of the theorem about the general section through $0$ and $\Gamma$, are satisfied. Let $W_{C,\Sigma}$ be the 
union of the images of the maps as in $(9)$. Then this is closed and by~\cite[Proposition 11.4]{Ko-Mo92} 
satisfies the requirements of the theorem.
\end{proof}

As an application we will consider the case of neighborhoods $Y \stackrel{f}{\la} X$, where $X$ is $cA_1$ or $cA_2$.

\subsection{The $cA_1$ case.}

This case behaves especially well as shown by the next theorem. 

\begin{theorem}\label{A1}
Let $C \subset Y \stackrel{f}{\la} X \ni 0$ be a divisorial extremal neighborhood with smooth center, such that $0 \in X$ is $cA_1$ type 
singularity. Let $0 \in T$ be the general section of $X$ through $0$, and $S$ the general section through the center $\Gamma$. Then
\begin{enumerate}
\item $T_Y=f^{\ast}T$ is a normal \textit{klt} surface with exactly one high index point.
\item $f$ can be obtained from the following diagram 
\[
\xymatrix{ 
            & Z \ar[dl]_{\nu} \ar[dr]^{\pi}\\
W \ar[dr]_g &   & Y \ar[dl]^f \\
            & X & }
\]
where all spaces involved are as in theorem~\ref{normal2}.
\item $S_Y \cong S$. In particular no semistable extremal neighborhood $f$ exists with $X$ $cA_1$.
\end{enumerate}
\end{theorem}
\begin{proof}
Parts $(2)$ and $(3)$ will follow from theorems~\ref{normal2},~\ref{nonnormal2} if we show that no semistable neighborhood exists with $X$ 
$cA_1$. Suppose there exists. Then by lemma~\ref{Pell}, such neighborhoods are in one to one correspondence with the solutions of the Pell 
equation $x^2+4y^2=4$, and the correspondence is given by $x=2n-n^{\prime}$, $y=n^{\prime}/2$, where $n,\; n^{\prime}$ are the indices of the 
two singular points of $Y$. But the only possibilities are $x=0, \; y=1$ and $x=2,\; y=0$. They give $n=1,\; n^{\prime}=2$, and 
$n=1,\; n^{\prime}=0$. Both are impossible and hence there is no semistable neighborhood of this type.

The only that is left is to show that $T_Y$ is normal. Suppose not. Then we are in the situation of theorems~\ref{nonnormal1},~\ref{nonnormal2}. 
Adopt the notation of these two theorems. Since $T$ is an $A_1$ DuVal singularity, by theorem~\ref{nonnormal2}.3, \[
\frac{m^2n}{(m,b+1)}(C_1-C_2)^2=2(\frac{mn}{(m,b+1)},\frac{m^2n}{(m,b+1)}(C_1-C_2)^2) \]
Therefore, there is a $\delta \in \mathbb{Z}$, such that \[
\delta \cdot \frac{1}{2} \frac{m^2n}{(m,b+1)}(C_1-C_2)^2=\frac{mn}{(m,b+1)} \]
and hence 
\begin{equation}
(C_1-C_2)^2=\frac{2}{m\delta}.
\end{equation}
Moreover, by theorem~\ref{nonnormal1}.3, \[
(C_1-C_2)^2=n^2[(C_1 \cdot C_2)^2-C_1^2 C_2^2].\]
Taking into consideration that $C_1 \cdot C_2 =1/m$, and that $C_1^2=a/n-b/m$,  we find that \[
\frac{2}{m\delta}=(C_1-C_2)^2=\frac{(n-ma-nb)^2}{m(m+amn-n^2b)}\]
and hence \[
\delta=\frac{2m+2amn-2n^2b}{(n-am-bn)^2} \in \mathbb{Z}.\]
Now \[
|\delta|\leq \frac{2m+2amn+2n^2b}{(n-am-bn)^2} \]
and it is easy to see that the above quantity is less than 1 if $b \geq 4$, and hence $\delta$ cannot be an integer. 
Therefore we must check only the cases $b=1,\; 2, \; 3$. 

I will only do the case $b=1$. The others are treated similarly.

If $b=1$, then $Q \in \overline{Z}\cong 1/m(1,-1)=1/m(1,m-1)$, and hence it is an $A_{m-1}$ DuVal singularity. 

Compute now $C_1^2$ and $C_2^2$. Suppose that in the minimal resolution of $\overline{Z}$, ${C_2^{\prime}}^2=-d$. Then 
\begin{gather*}
C_1^2=\frac{a}{n}-\frac{1}{m} \\
C_2^2=-d+\frac{n-a}{n}+\frac{m-1}{m}=2-d-\frac{a}{n}-\frac{1}{m} 
\end{gather*}
Moreover, $C_1 \cdot C_2 =1/m$, and hence \[
(C_1-C_2)^2=-\frac{4}{m}+2-d \]
But now from $(10)$ it follows that \[
\frac{2}{m(C_1-C_2)^2}=\frac{2}{-4+2m-md} \in \mathbb{Z} \]
It is now easy to see that the only possibilities are $(m,d)=(5,1)$, or $(3,1)$, or $(6,1)$, or $(2,1)$. In any case $d=1$. 
Hence the part of the extended dual graph of $\overline{Z}$ that corresponds to $C_1$, $C_2$ and $Q$ is \[
\stackrel{-1}{\underset{C_1}{\bullet}}\mbox{\noindent ---}\stackrel{-2}{\underset{E_1}{\circ}}\mbox{\noindent ---} \cdots 
\mbox{\noindent ---}\stackrel{-2}{\underset{E_m}{\circ}}\mbox{\noindent ---}\stackrel{-1}{\underset{C_2}{\bullet}} \]
This configuration should be contractible. But \[
(C_1+\sum_{i=1}^m E_i +C_2)^2=-1-2m-1+2+2(m-1)+2=0\]
which is impossible since it must be negative. Hence $T_Y$ must be normal and theorem~\ref{A1} follows.
\end{proof}

\subsection{The $cA_2$ case}

We will now classify semistable extremal neighborhoods $Y \stackrel{f}{\la} X$, such that $X$ is $cA_2$ type. 
Again by using corollary~\ref{Pell}, we see that such neighborhoods are in one to one correspondence with the integer 
solutions of the Pell equation $x^2+3y^2=4$, and the correspondence is given by $x=2n-n^{\prime}$, $y=n^{\prime}/3$. 
The only solution is $x=1,\; y=1$, which gives $n=2$ and $n^{\prime}=3$. Moreover, $dd^{\prime}=3$. 
Hence from this and theorem~\ref{semistable1}.5, it follows that $d^{\prime}=1$ and $d=3$. From theorem~\ref{semistable1}.4 
it also follows that $3a+2a^{\prime}=7$
and hence $a=1$ and $a^{\prime}=2$. Therefore, $(P\in Z)\cong 1/12(1,5) \cong (xy-z^6=0)/\mathbb{Z}_2(1,-1,1)$, 
and  $(P^{\prime} \in Z)\cong 1/9(1,5) \cong (xy-z^3=0)/\mathbb{Z}_3(1,-1,2)$.

Now let $S$ be the general section of $X$ through $\Gamma$. We would like to find conditions under which a semistable 
neighborhood with these data does not exist. Suppose that $S$ is $A_m$ type DuVal singularity and that $\Gamma$ 
intersects the $s$ part of the dual graph of $S$. Then by theorem~\ref{semistable2}, it follows that 
\begin{gather*}
2k+3k^{\prime}=m+1 \\
k+k^{\prime}=s
\end{gather*}
They give that $3s-m-1=k$, and hence if it happens that $s\leq (m+1)/3$, such a $k$ does not exist and hence no semistable neighborhood 
as well. 

Hence we have shown the following.
 
\begin{proposition}
Semistable divisorial extremal neighborhoods $C \subset Y \stackrel{f}{\la} X \ni 0$, such that $0 \in X$ is a $cA_2$ cDV 
singularity, are in one to one correspondence with $\mathbb{Q}$-Gorenstein smoothings of the surface germ $C \subset Z$, 
such that the singular locus of $Z$ is just two points $(P\in Z)\cong (xy-z^6=0)/\mathbb{Z}_2(1,-1,1)$, and 
$(P^{\prime} \in Z)\cong (xy-z^3=0)/\mathbb{Z}_3(1,-1,2)$. 

Let $S$ be the general section of $X$ through $\Gamma$. $S$ is an $A_m$ DuVal singularity for some $m$, and suppose that $\Gamma$ 
intersects the $s$ part of the dual graph of $S$. Then if $s \leq (m+1)/3$, no semistable contraction with this data exists.
\end{proposition}